\documentclass{article}
\usepackage{wrapfig}
\usepackage{amsmath}
\usepackage{amssymb}
\usepackage{fullpage}
\usepackage{epsfig}

\begin{document}
\title{A Variational Surface-Evolution Perspective for Optimal Transport
between Densities with Differing Compact Support}
\author{Anthony Yezzi}
\date{}
\maketitle

\begin{abstract}
We examine the optimal mass transport problem in $\mathbb{R}^{n}$
between densities having independent compact support by considering
the geometry of a continuous interpolating support boundary in space-time
within which the mass density evolves according to the fluid dynamical
framework of Benamou and Brenier. We treat the geometry of this space--time
embedding in terms of points, vectors, and sets in $\mathbb{R}^{n+1}\!=\mathbb{R}\times\mathbb{R}^{n}$
and blend the mass density and velocity as well into a space-time
solenoidal vector field ${\bf W}\;|\;{\bf \Omega\subset}\mathbb{R}^{n+1}\!\to\mathbb{R}^{n+1}$
over compact sets ${\bf \Omega}$ . We then formulate a coupled gradient
descent approach containing separate evolution steps for $\partial{\bf \Omega}$
and ${\bf W}$.
\end{abstract}

\section{Introduction}

Optimal mass transport (OMT) has a long history, beginning with Gaspard
Monge \cite{GM} in 1781, and put into a more modern form solvable via linear
programming by Leonid Kantorovich \cite{LK,GangboMcCann,villani1,villani2,Evans1999,Rachev1998}.
In recent years, OMT has undergone a huge surge, with many diverse
applications including signal/image processing, computer vision, machine
learning, data analysis, meteorology, statistical physics, quantum
mechanics, and network theory \cite{Arjovsky2017,CarMaa14,Haker2004,MitMie16,rr,Statement2020}. 

Our interest in the present work is extending the Benamou and Brenier
computational fluid dynamical (CFD) to OMT where the domains have
compact support. We recall that in the seminal work \cite{BB}, Benamou
and Brenier compute the Wasserstein 2-metric $W_{2}$ via the minimization
of kinetic energy subject to a continuity constraint. They moreover
compute the optimal path (i.e., geodesic) in the space of probability
densities \cite{Otto}. In image processing, this gives a natural
interpolation path between two images, where the intensity is treated
as a ``generalized mass.'' This is important for problems in deformable
registration and image warping; see \cite{omtEx3,haker} and the references
therein..

Now the method of Benamou and Brenier does not take into account that
the regions of interest in images are in fact compact, and so the
computations cannot be performed on a regular grid and crucially boundary
conditions must be imposed. It is this issue that motivated the present
work. To treat the problem of OMT on densities with compact support,
we propose a synergy of methods from OMT as well as level set evolutions
\cite{osher,sethian}. Level set methods are a powerful way of implementing
interface (boundary) evolution problems and thus have become a standard
for a number of approaches in computer vision for segmentation, so-called
\emph{active contour} methods; see \cite{sapiro} and the references
therein. 

This synergy of level sets and OMT, we believe to be novel, and to
potentially have a number of applications, in particular to digital
pathology.

\section{Spatiotemporal hypersurface}

We use \textbf{bold notation} exclusively to denote space-time points,
vectors, and sets and non-bold notation for similar entities in time
or space only, as well as for scalar variables. Accordingly, ${\bf X}$
will represent an arbitrary spatiotemporal point in $\mathbb{R}^{n+1}$
\[
{\bf X}=(t,x)=(\underbrace{X_{0}}_{t},\underbrace{X_{1},X_{2},\ldots,X_{n}}_{x})
\]
which pairs a temporal coordinate $t\in\mathbb{R}$ with a spatial
coordinate $x\in\mathbb{R}^{n}$. We will denote the spatiotemporal
basis vectors by ${\bf e_{0}},{\bf e}_{1},{\bf e}_{2},\ldots,{\bf e}_{n}$
so that
\[
{\bf X}=\underbrace{X_{0}\,{\bf e}_{0}}_{t}+\underbrace{X_{1}\,{\bf e}_{1}+X_{2}\,{\bf e}_{2}+\cdots+X_{n}\,{\bf e}_{n}}_{x}
\]

\subsection{Assumptions}

\paragraph{Compact support}

We consider a spatiotemporal density $\rho$ in the form of a positive
scalar function ${\bf X}\to\rho({\bf X})$ in $\mathbb{R}^{n+1}$
whose restriction to the $t\!=\!0$ hyperplane matches a given initial
spatial density $\rho_{0}$ and whose restriction to the $t\!=\!1$
hyperplane matches a given spatial target density $\rho_{1}$. We
assume the initial spatial density $\rho_{0}$ has compact support
$\Omega_{0}\subset\mathbb{R}^{n}$ and that the target spatial density
$\rho_{1}$ has compact support $\Omega_{1}\subset\mathbb{R}^{n}$.
We further assume that the full spatiotemporal density $\rho$ also
has compact support ${\bf \Omega}\subset\mathbb{R}^{n+1}$ in space-time,
sandwiched between the hyperplanes $t\!=\!0$ and $t\!=\!1$, which
may be constructed by a continuous family of intermediate compact
spatial supports $\Omega_{[t]}\subset\mathbb{R}^{n}$ with $\Omega_{[0]}=\Omega_{0}$
and $\Omega_{[1]}=\Omega_{1}$ as follows.
\begin{equation}
{\bf \Omega}=\left\{ {\bf X}=(t,x)\;|\;0\le t\le1\,,\,x\in\Omega_{[t]}\right\} \label{eq:compact}
\end{equation}

\paragraph{Balanced density}

We assume the initial and target spatial densities $\rho_{0}$ and
$\rho_{1}$ both have unit mass. We impose a similar constraint on
the spatiotemporal density $\rho$, summarizing these assumptions
as follows.\footnote{The reason we don't start with the stronger constraint $\int_{\Omega_{[t]}}\rho(x,t)\,dx=1$
for all $0\le t\le1$ is that both this as well as the weaker total
spatiotemporal mass constraint are global mass preservation constraints
that will automatically be satisfied later when we impose the much
stronger local mass preservation constraint. The main point in even
presenting the total spatiotemporal constraint here is to reinforce
the embedded space-time interpretation of the problem, thereby allowing
us to interpret $\rho$ as a unit spatiotemporal mass density directly
in $\mathbb{R}^{n+1}$.}
\begin{equation}
\int_{\Omega_{0}}\rho_{0}(x)\,dx=\int_{\Omega_{1}}\rho_{1}(x)\,dx=\int_{{\bf \Omega}}\rho({\bf X})\,d{\bf X}=1\label{eq:balanced}
\end{equation}

\paragraph{Smoothness}

We assume the initial and target spatial densities $\rho_{0}$ and
$\rho_{1}$ are differentiable \emph{within their support}\footnote{We do not require $\rho_{0}$ and $\rho_{1}$ to be zero along the
boundaries of their supports, as such they may discontinuously drop
to zero across the spatial boundaries $\Gamma_{0}$ and $\Gamma_{1}$
respectively.} $\Omega_{0}$ and $\Omega_{1}$ and that the spatiotemporal density
$\rho$ is differentiable \emph{within its support}\footnote{Nor do we require $\rho$ to be zero along the boundary of its support
(a necessary freedom along the flat temporal faces ${\bf \Gamma}_{0}$
and ${\bf \Gamma}_{1}$ , else we could not impose $\rho=\rho_{o}$
and $\rho=\rho_{1}$ along these components).} ${\bf \Omega}$. We also assume that the portion of the spatiotemporal
boundary $\partial{\bf \Omega}$ that lies strictly within $0\!<\!t\!<\!1$,
which we denote by ${\bf \Gamma}$, is differentiable. The remaining
portions of $\partial{\bf \Omega}$ are provided by the embeddings
of $\Omega{}_{0}$ and $\Omega_{1}$ within the hyperplanes $t\!=\!0$
and $t\!=\!1$ to form two flat \emph{temporal faces} of ${\bf \Omega}$,
which we denote by ${\bf \Gamma}_{0}$ and ${\bf \Gamma}_{1}$. 

\paragraph{Piece-wise smooth boundary }

As such $\partial{\bf \Omega={\bf \Gamma}}\cup{\bf \Gamma}_{0}\cup{\bf \Gamma}_{1}$
will have the form of a compact hypersurface in $\mathbb{R}^{n+1}$
with a well defined outward normal everywhere except along the borders
of the two temporal faces where ${\bf \Gamma}_{0}$ and ${\bf \Gamma}_{1}$
connect to the intervening surface ${\bf \Gamma}$. We may also describe
the intervening spatiotemporal boundary component ${\bf \Gamma}$
as the \emph{swept out surface} generated by embedding the boundaries
$\Gamma_{[t]}=\partial\Omega_{[t]}$ of the deforming spatial supports
$\Omega_{[t]}$ into the corresponding temporal hyperplanes. The spatiotemporal
boundary notation and decomposition is summarized as follows.
\begin{equation}
\partial{\bf {\bf \Omega}}=\underbrace{\left\{ {\bf X}=(0,x)\;|\;x\in\Omega_{0}\right\} }_{\underset{\mbox{(temporal face)}}{{\bf \Gamma}_{0}}}\cup\underbrace{\left\{ {\bf X}=(1,x)\;|\;x\in\Omega_{1}\right\} }_{\underset{\mbox{(temporal face)}}{{\bf \Gamma}_{1}}}\cup\underbrace{\left\{ {\bf X}=(t,x)\;|\;0\le t\le1\,,\,x\in\overbrace{\partial\Omega_{[t]}}^{\Gamma_{[t]}}\right\} }_{\underset{\mbox{(swept out surface)}}{{\bf \Gamma}}}\label{eq:boundary}
\end{equation}

\subsection{Local geometry}

In this section we explore the relationship between the local geometry
of the spatial support boundary $\Gamma_{[t]}$ and the swept out
spatiotemporal boundary ${\bf \Gamma}$.

\paragraph{Parameterization}

Let $s=(s_{1},\ldots,s_{n-1})$ represent isothermal coordinates with
unit speed at a point $x\in$$\Gamma_{[t]}$ along the boundary of
the spatial support $\Omega_{[t]}\subset\mathbb{R}^{n}$ at time $t$
strictly between 0 and 1, The Riemannian metric tensor of $\Gamma_{[t]}$
in these coordinates is therefore the $(n-1)\times(n-1)$ identity
matrix at the point $x$ with $n-1$ orthonormal tangent vectors $\frac{\partial\Gamma_{[t]}}{\partial s_{k}}\in\mathbb{R}^{n}$
for $k=1,\ldots,n-1$. From (\ref{eq:boundary}) we know that the
corresponding spatiotemporal point ${\bf X}=(t,x)$ belongs to portion
${\bf \Gamma}$ of the spatiotemporal support boundary $\partial{\bf \Omega}$,
which may be locally parameterized as follows.
\begin{equation}
{\bf \Gamma}(t,s)=\left(t,\Gamma_{[t]}(s)\right)=(\underbrace{\Gamma_{0}}_{t},\underbrace{\Gamma_{1},\Gamma_{2},\ldots,\Gamma_{n}}_{\Gamma_{[t]}(s)})\label{eq:parameterization}
\end{equation}

\paragraph{Unit normal}

In these coordinates, we compute the following $n$ tangent vectors
to ${\bf \Gamma}$ in $\mathbb{R}^{n+1}$
\begin{equation}
\frac{\partial{\bf \Gamma}}{\partial t}=\left(1,\frac{\partial\Gamma_{[t]}}{\partial t}\right)\quad\mbox{and}\quad\frac{\partial{\bf \Gamma}}{\partial s_{k}}=\left(0,\frac{\partial\Gamma_{[t]}}{\partial s_{k}}\right),\text{\ensuremath{\quad k=1,\ldots,n-1}}\label{eq:tangents}
\end{equation}
Since the unit outward normal $N\in\mathbb{R}^{n}$ to the spatial
boundary $\Gamma_{[t]}$ is orthogonal to the $n-1$ spatial tangent
vectors $\frac{\partial\Gamma_{[t]}}{\partial s_{k}}\in\mathbb{R}^{n}$,
it is clear to see that $(\alpha,N)\in\mathbb{R}^{n+1}$ will be orthogonal
to the $n-1$ spatiotemporal tangent vectors $\frac{\partial{\bf \Gamma}}{\partial s_{k}}$
expressed in (\ref{eq:tangents}) for any choice of scalar $\alpha$.
Orthogonality to the additional spatiotemporal tangent vector $\frac{\partial{\bf \Gamma}}{\partial t}$
expressed in (\ref{eq:tangents}) as well requires $\alpha=-\frac{\partial\Gamma_{[t]}}{\partial t}\cdot N$,
yielding the following construction of the spatiotemporal unit outward
normal ${\bf N}\in\mathbb{R}^{n+1}$.

\begin{equation}
\text{{\bf N}}=\begin{cases}
\frac{\left(-\frac{\partial\Gamma_{[t]}}{\partial t}\cdot N,N\right)}{\sqrt{1+\left(\frac{\partial\Gamma_{[t]}}{\partial t}\cdot N\right)^{2}}}, & \text{{\bf X}}\in{\bf \Gamma}\\
+{\bf e}_{0}=+(1,0,\ldots,0), & \text{{\bf X}}\in{\bf \Gamma}_{1}\\
-{\bf e}_{0}=-(1,0,\ldots,0),, & \text{{\bf X}}\in{\bf \Gamma}_{0}
\end{cases}\label{eq:normal}
\end{equation}

\paragraph{Metric tensor}

Using the tangent vectors (\ref{eq:tangents}) we may express the
Riemannian metric tensor at the point ${\bf X}$ in the form of the
following $n\times n$ matrix
\begin{equation}
\underbrace{{\bf I}(s_{1},\ldots,s_{n-1},t)}_{\mbox{I-fundamental form}}=\begin{bmatrix}{\cal I} & \left(\frac{\partial\Gamma_{[t]}}{\partial s}\right)^{T}\frac{\partial\Gamma_{[t]}}{\partial t}\\
\frac{\partial\Gamma_{[t]}}{\partial t}^{T}\left(\frac{\partial\Gamma_{[t]}}{\partial s}\right) & 1+\frac{\partial\Gamma_{[t]}}{\partial t}\cdot\frac{\partial\Gamma_{[t]}}{\partial t}
\end{bmatrix}\label{eq:metric}
\end{equation}
where $\frac{\partial\Gamma_{[t]}}{\partial s}$ denotes the $n\times(n-1)$
matrix whose columns consist of the orthonormal tangent vectors to
the spatial boundary $\Gamma_{[t]}$ at the point $x$. Using the
determinant formula, $\det\begin{bmatrix}A & u\\
v^{T} & \alpha
\end{bmatrix}=\alpha\det A-v^{T}\mbox{adj}A\,u$ (for any matrix $A$, vector $u$ and $v$, and scalar $\alpha$),
we may compute.
\[
\det{\bf I}=\left(1+\frac{\partial\Gamma_{[t]}}{\partial t}\cdot\frac{\partial\Gamma_{[t]}}{\partial t}\right)-\underbrace{\sum_{k=1}^{n-1}\left(\frac{\partial\Gamma_{[t]}}{\partial s_{k}}\cdot\frac{\partial\Gamma_{[t]}}{\partial t}\right)^{2}}_{\left\Vert \frac{\partial\Gamma_{[t]}}{\partial t}\right\Vert ^{2}-\left(\frac{\partial\Gamma_{[t]}}{\partial t}\cdot N\right)^{2}}=1+\left(\frac{\partial\Gamma_{[t]}}{\partial t}\cdot N\right)^{2}
\]

\paragraph{Area element}

The relationship between the area element $dS_{[t]}$ of the spatial
boundary surface $\Gamma_{[t]}$, the temporal element $dt$, and
the area element $d{\bf S}$ of the swept out spatiotemporal hypersurface
${\bf \Gamma}$ can be expressed via the square root of the determinant
of the first fundamental form shown above.
\begin{equation}
d{\bf S}=\sqrt{1+\left(\frac{\partial\Gamma_{[t]}}{\partial t}\cdot N\right)^{2}}dS_{[t]}\,dt\label{eq:area}
\end{equation}

\paragraph{Normal variations}

Finally, using the parameterization (\ref{eq:parameterization}) for
a variation $\delta{\bf {\bf \Gamma}}$ allows us to relate a variation
of the swept out spatiotemporal hypersurface to time parameterized
variations $\delta\Gamma_{[t]}$ of the spatial support boundaries
as follows $\delta{\bf {\bf \Gamma}}=\left(0,\delta\Gamma_{[t]}\right)$.
Combining this with (\ref{eq:normal}) yields the following relationship
between variations of $\Gamma_{[t]}$ and variations of ${\bf \Gamma}$
along their respective normal directions.

\[
\delta{\bf {\bf \Gamma}}\cdot{\bf N}=\frac{\delta\Gamma_{[t]}\cdot N}{\sqrt{1+\left(\frac{\partial\Gamma_{[t]}}{\partial t}\cdot N\right)^{2}}}
\]
If we now further combine this with (\ref{eq:area}) we see that the
normal variation of the spatiotemporal hypersurface measured against
its spatiotemporal area element matches the normal variation of the
corresponding spatial boundary measured by its respective area element
and time element.
\begin{equation}
\left(\delta{\bf {\bf \Gamma}}\cdot{\bf N}\right)\,d{\bf S}=\left(\delta\Gamma_{[t]}\cdot N\right)\,dS_{[t]}\,dt\label{eq:normal-measures}
\end{equation}

\section{Spatiotemporal formulation of optimal mass transport}

The fluid dynamical framework of Benamou and Brenier, considers two
time evolving entities, a scalar mass density $\rho(t,x)$ and a velocity
field $v(t,x)\in\mathbb{R}^{n}$ which are coupled by the local mass
preserving continuity constraint $\frac{\partial\rho}{\partial t}+\nabla_{x}\cdot(\rho v)=0$. 

Notice that in the standard manner regarding $\rho v$ as an ordinary
3-vector, the continuity constraint means that density and spatial
momentum form a 4-vector with respect to the standard Minkowski metric
in the standard physics setting \cite{susskind}. We will exploit
this observation and develop an equivalent formulation using a single
space-time solenoidal vector field ${\bf W}$ with simple normal boundary
conditions along the spatiotemporal hypersurface $\partial{\bf \Omega}$
that enable convenient numerical solutions of PDE's directly within
on an $n+1$ dimensional space-time grid, with no need to treat the
temporal and spatial dimensions separately or differently. 

\subsection{Spatiotemporal advection field ${\bf U}$}

We begin by noting that in the combined spatiotemporal variable ${\bf X}=(t,x)$,
the continuity constraint $\frac{\partial\rho}{\partial t}+\nabla_{x}\cdot(\rho v)=0$
can be written as
\begin{equation}
\nabla\rho\cdot{\bf U}+\rho\nabla\cdot{\bf U}=0\label{eq:transport}
\end{equation}
where $\nabla$ and $\nabla\cdot$ represent the full spatiotemporal
gradient and divergence operators in $\mathbb{R}^{n+1}$ and where
${\bf U}\in\mathbb{R}^{n+1}$ denotes the following vector field.
\[
{\bf U}({\bf X})\dot{=}\left(1,v\right)=(\underbrace{U_{0}}_{1},\underbrace{U_{1},U_{2},\ldots,U_{n}}_{v})
\]
Note that ${\bf U}$ is tangent to the characteristics of this linear
first order PDE (\ref{eq:transport}) in $\rho$ and therefore defines
the trajectories along which mass is transported across space-time.
Since the spatiotemporal hypersurface $\partial{\bf \Omega}$, more
specifically its swept out portion ${\bf \Gamma}$, defines the boundary
of the evolving support for $\rho$, we can conclude that these advection
trajectories must flow along the hypersurface ${\bf \Gamma}$ itself,
never across it (neither inward nor outward). In other words, mass
cannot be transported outside of its support neither forward in time
(which excludes outward-flowing characteristics) nor backward in time
(which excludes inward-flowing characteristics). This leads to the
boundary condition ${\bf U}\cdot{\bf N}=0$ along ${\bf \Gamma}$,
which may be written in terms of $v$, $N$, and $\frac{\partial\Gamma_{[t]}}{\partial t}$
as follows. 
\begin{equation}
\underbrace{v\cdot N=\frac{\partial\Gamma_{[t]}}{\partial t}\cdot N}_{({\bf U}\cdot{\bf N}=0)}\label{eq:support-evolution}
\end{equation}
If we plug this constraint between the support evolution $\frac{\partial\Gamma_{[t]}}{\partial t}$
and the velocity field $v$ into (\ref{eq:normal}) we obtain the
following alternative expression for the outward unit normal ${\bf N}\in\mathbb{R}^{n+1}$
of the swept out hypersurface ${\bf \Gamma}$ in terms of the outward
normal $N\in\mathbb{R}^{n}$ of the support boundary $\Gamma_{[t]}$
in space at time $t$.
\begin{align}
{\bf N}({\bf X}) & =\frac{\left(-v\cdot N,N\right)}{\sqrt{1+\left(v\cdot N\right)^{2}}},\quad{\bf X}\in{\bf \Gamma}\label{eq:normal-v}
\end{align}

\subsection{Solenoidal vector field ${\bf W}$}

While the use of the advection field ${\bf U}$ merges the spatial
and temporal derivatives into a single derivative $\mathbb{R}^{n+1}$
operator in (\ref{eq:support-evolution}), it still keeps the density
variable $\rho$ separate. We now merge these two entities by defining
another spatiotemporal vector field whose first (temporal) component
represents the mass density $\rho({\bf X})$, and whose remaining
components represent the momentum vector ${\rm p}({\bf X})=\rho v$
in $\mathbb{R}^{n}.$
\begin{equation}
{\bf W}({\bf X})=\rho{\bf U}=\left(\,\rho\,,\,{\rm p}\,\right)=(\underbrace{W_{0}}_{\rho},\underbrace{W_{1},W_{2},\ldots,W_{n}}_{\mbox{momentum }{\rm p}=\rho v})\label{eq:merged}
\end{equation}
However, rather than considering (\ref{eq:merged}) to be the definition
of ${\bf W}$ in terms of the density $\rho$ and momentum ${\rm p}$,
we instead consider it in reverse to be the definition of $\rho$
and ${\rm p}$ in terms of the space-time vector field ${\bf W}$
subject to the continuity constraint (\ref{eq:transport}) which now
simplifies to a coordinate-free solenoidal condition on ${\bf W}$.
\begin{equation}
\nabla\cdot{\bf W}=0\label{eq:solenoidal}
\end{equation}
Multiplying the boundary condition ${\bf U}\cdot{\bf N}=0$ along
${\bf \Gamma}$ presented in (\ref{eq:support-evolution}) by $\rho$
yields a similar vanishing flux condition for ${\bf W}$ across the
swept out portion of the spatiotemporal hypersurface $\partial{\bf \Omega}$.
For the remainder of $\partial{\bf \Omega}$, we combine ${\bf N}=\text{\ensuremath{\pm}}{\bf e}_{0}$
from (\ref{eq:normal}) with (\ref{eq:merged}) to obtain flux conditions
for ${\bf W}$ along the temporal faces ${\bf \Gamma}_{0}$ and ${\bf \Gamma}_{1}$
as well, in terms of the known starting and target densities $\rho_{0}$
and $\rho_{1}$. 

The combination of these constraints is easily summarized now in terms
of ${\bf W}$ and its spatiotemporal domain ${\bf \Omega}$. Namely
we seek a \emph{solenoidal vector field} ${\bf W}$ within ${\bf \Omega}$
with the following prescribed flux conditions along the full spatiotemporal
domain boundary $\partial{\bf {\bf \Omega}}$.
\begin{equation}
{\bf W}\cdot{\bf N}=\begin{cases}
\;\;0\;\;, & {\bf X}\in{\bf \Gamma}\\
-\rho_{0}, & {\bf X}\in{\bf \Gamma}_{0}\\
+\rho_{1}, & {\bf X}\in{\bf \Gamma}_{1}
\end{cases}\label{eq:flux}
\end{equation}

\subsection{Extended velocity ${\bf V}$}

\paragraph{Local kinetic energy}

Before setting up the variational problem we seek an expression for
the local measure of kinetic energy ${\rm T}(t,x)$
\[
{\rm T}({\bf X})=\frac{1}{2}\rho\,v\cdot v
\]
in terms of the solenoidal field ${\bf W}$. If we express this purely
in terms of ${\bf W}$, we obtain the following expression which,
unfortunately, is not coordinate free.

\[
{\rm T}({\bf X})=\frac{1}{2}\left(\frac{\overbrace{{\bf W}\cdot{\bf W}}^{\rho^{2}(1+v\cdot v)}}{\underbrace{{\bf W}\cdot{\bf e}_{0}}_{\rho}}-\underbrace{{\bf W}\cdot{\bf e}_{0}}_{\rho}\right)
\]

\paragraph{Extended velocity}

We may resolve this by introducing the following \emph{extended velocity}
field ${\bf V}\in\mathbb{R}^{n+1}$ which extends the transport velocity
$v$ from $\mathbb{R}^{n}$ into $\mathbb{R}^{n+1}$ by adding a temporal
component equal to $-\frac{1}{2}\|v\|^{2}$ as follows
\begin{equation}
{\bf V}\doteq\left(-\frac{1}{2}\|v\|^{2}\,,\,v\right)=(\underbrace{V_{0}}_{\frac{-\|v\|^{2}}{2}},\underbrace{V_{1},V_{2},\ldots,V_{n}}_{v})\label{eq:V}
\end{equation}
Notice that, just like the advection field ${\bf U}$, the extended
velocity ${\bf V}$ depends only upon the spatial velocity $v$ itself,
and therefore contains no additional information. We may now express
the local kinetic energy compactly and coordinate free in terms of
the solenoidal field ${\bf W}$ and the extended velocity field ${\bf V}$
by their inner product.
\[
{\rm T}={\bf W}\cdot{\bf V}
\]
We summarize our notation for the three spatial-temporal vector fields
in $\mathbb{R}^{n+1}$ as follows.

\smallskip{}

\[
\begin{array}{ccc}
\mbox{advection field} & \mbox{solenoidal field} & \mbox{extended velocity}\\
{\bf U}=(1,v) & {\bf W}=\rho{\bf U}=(\rho,{\rm p}) & {\bf V}=\left(-\frac{1}{2}\|v\|^{2}\,,\,v\right)
\end{array}
\]

\subsection{Variational formulation }

\subsubsection{Action integral}

We begin by constructing the action integral to be minimized over
${\bf \Omega}\subset\mathbb{R}^{n+1}$ in terms of the solenoidal
field ${\bf W}$ as follows. 
\[
\int_{0}^{1}\int_{\Omega_{[t]}}\overbrace{\frac{1}{2}\rho\|v\|^{2}}^{T}\,dx\,dt=\int_{{\bf \Omega}}{\bf W}\cdot{\bf V}\,d{\bf X}
\]
Note that the two unknowns are only ${\bf W}$ and its support ${\bf \Omega}$
(or equivalently the swept out boundary ${\bf \Gamma}$) even though
we have expressed the action compactly also in terms of ${\bf V}$.
We may compute ${\bf V}$ directly from ${\bf W}$ 
\begin{equation}
{\bf V}\,\doteq\left(-\frac{1}{2}\|v\|^{2}\,,\,v\right)={\bf U}-\frac{1}{2}\left({\bf U}\cdot{\bf U}+1\right){\bf e}_{0}=\frac{{\bf W}}{{\bf W}\cdot{\bf e}_{0}}-\frac{1}{2}\left(\frac{{\bf W}\cdot{\bf W}}{\left({\bf W}\cdot{\bf e}_{0}\right)^{2}}+1\right){\bf e}_{0}\label{eq:V_W}
\end{equation}
with the following compatible flux conditions obtained by plugging
in (\ref{eq:flux}). 
\begin{align*}
{\bf V}\cdot{\bf N} & =\begin{cases}
-\left(1+\frac{1}{2}\|v\|^{2}\right)\left({\bf N}\cdot{\bf e}_{0}\right), & {\bf X}\in{\bf \Gamma}\\
+\frac{1}{2}\|v\|^{2}, & {\bf X}\in{\bf \Gamma}_{0}\\
-\frac{1}{2}\|v\|^{2}, & {\bf X}\in{\bf \Gamma}_{1}
\end{cases}
\end{align*}
Incorporating the solenoidal (mass preservation) constraint (\ref{eq:solenoidal})
through a Lagrange multiplier $\lambda$ over ${\bf {\bf \Omega}}$
and the flux constraints through an additional Lagrange multiplier
$\alpha$ along the boundaries $\partial{\bf \Omega}={\bf \Gamma}_{0}\cup{\bf \Gamma}_{1}\cup{\bf \Gamma}$.
\begin{equation}
E=\int_{{\bf \Omega}}{\bf W}\cdot{\bf V}\,+\underbrace{\lambda\,\nabla\cdot{\bf W}}_{\stackrel{\mbox{solenoidal}}{\mbox{constraint}}}\,d{\bf X}+\underbrace{\int_{{\bf \Gamma}}\alpha{\bf W}\cdot{\bf N}\,dS+\int_{{\bf \Gamma}_{0}}\alpha\left(\overbrace{{\bf W}\cdot{\bf N}}^{-{\bf W}\cdot{\bf {\bf e}_{0}}}+\rho_{0}\right)\,dS+\int_{{\bf \Gamma}_{1}}\alpha\left(\overbrace{{\bf W}\cdot{\bf N}}^{{\bf W}\cdot{\bf e}_{0}}-\rho_{1}\right)\,dS}_{\mbox{flux constraints}}\label{eq:energy}
\end{equation}

\subsubsection{First variation}

Next, we compute the variation of (\ref{eq:energy}) to obtain (see
Appendix \ref{app:variation})
\begin{align}
\delta E= & \int_{{\bf \Omega}}\delta{\bf W}\cdot\left({\bf V}-\nabla\lambda\right)\;d{\bf X}+\int_{{\bf \Gamma}}{\bf W}\cdot\left({\bf V}+\nabla_{S}\alpha\right)\delta{\bf \Gamma}\cdot{\bf N}\,dS\label{eq:variation}\\
 & +\underbrace{\int_{{\bf \Gamma}}\left(\alpha+\lambda\right)\,\delta{\bf W}\cdot{\bf N}\,dS-\int_{{\bf \Gamma}_{0}}\left(\alpha+\lambda\right)\,\delta{\bf W}\cdot{\bf e}_{0}\,dS+\int_{{\bf \Gamma}_{1}}\left(\alpha+\lambda\right)\,\delta{\bf W}\cdot{\bf e}_{0}\,dS}_{\int_{\partial{\bf \Omega}}\left(\alpha+\lambda\right)\,\delta{\bf W}\cdot{\bf N}\,dS}\nonumber 
\end{align}

\paragraph{Necessary optimality condition}

Optimality for ${\bf W}$ can only be achieved if we can solve
\begin{equation}
\nabla\lambda={\bf V}\label{eq:hj_compact}
\end{equation}
for $\lambda$ within the interior of ${\bf \Omega}$ (additional
conditions are required along the boundary to annihilate the boundary
integrals as well). If we separately equate the spatial and temporal
components 

\[
\overbrace{\Bigl(\underbrace{\frac{\partial\lambda}{\partial X_{0}}}_{\frac{\partial\lambda}{\partial t}},\underbrace{\frac{\partial\lambda}{\partial X_{1}},\frac{\partial\lambda}{\partial X_{2}},\ldots,\frac{\partial\lambda}{\partial X_{n}}}_{\frac{\partial\lambda}{\partial x}}\Bigr)}^{\nabla\lambda}=\overbrace{\Bigl(\underbrace{V_{0}}_{\frac{-\|v\|^{2}}{2}},\underbrace{V_{1},V_{2},\ldots,V_{n}}_{v}\Bigr)}^{{\bf V}}
\]
then we see that (\ref{eq:hj_compact}) amounts to a more compact,
coordinate-free expression of the well known Hamilton Jacobi equation.
\begin{equation}
\frac{\partial\lambda}{\partial t}+\frac{1}{2}\left\Vert \frac{\partial\lambda}{\partial x}\right\Vert ^{2}=0\label{eq:hamilton-jacobi}
\end{equation}
However, (\ref{eq:hj_compact}) will not admit a solution unless ${\bf V}$
is a conservative (irrotational) vector field. As such, the gradient
must be related to the non-conservative (solenoidal) portion of the
extended velocity field ${\bf V}$. 

\paragraph{Helmoltz decomposition}

Accordingly, we consider the Helmholtz decomposition to express ${\bf V}$
as the sum of two vector fields
\begin{equation}
\text{{\bf V}}=\underbrace{\text{{\bf V}}^{\parallel}}_{\stackrel{\mbox{irrotational}}{(\mbox{curl-free)}}}+\underbrace{{\bf V}^{\perp}}_{\stackrel{\mbox{solenoidal}}{(\mbox{divergence-free})}}\label{eq:helmholtz}
\end{equation}
where ${\bf V}^{\perp}$ denotes the divergence free component ($\nabla\cdot{\bf V}^{\perp}=0$)
and where ${\bf V}^{\parallel}$ denotes the curl free component ($\nabla\times{\bf V}^{\parallel}=0$)
which can be written as the gradient of a scalar potential function
${\bf V}^{\|}=\nabla\phi$. However, this decomposition is not unique
over compact domains. We obtain the decomposition by solving the following
Poisson equation for a scalar potential function $\phi$.
\begin{align}
\Delta\phi & =\nabla\cdot{\bf V},\quad{\bf X}\in{\bf \Omega}\label{eq:poisson}
\end{align}
and can therefore parameterize the set of all possible decompositions
by the choice of imposed boundary conditions along $\partial{\bf \Omega}$.
To maintain the initial and final density constraints, Neumann boundary
conditions must be imposed on the two temporal faces ${\bf \Gamma}_{0}$
and ${\bf \Gamma}_{1}$

\begin{align}
\nabla\phi\cdot{\bf N} & =\overbrace{{\bf V}\cdot{\bf N}}^{-{\bf V}\cdot{\bf e}_{0}}=+\frac{1}{2}\|v\|^{2},\;{\bf X}\in{\bf \Gamma}_{0}\label{eq:temporal_bc}\\
\nabla\phi\cdot{\bf N} & =\underbrace{{\bf V}\cdot{\bf N}}_{+{\bf V}\cdot{\bf e}_{0}}=-\frac{1}{2}\|v\|^{2},\;{\bf X}\in{\bf \Gamma}_{1}\nonumber 
\end{align}
thereby leaving the decomposition (\ref{eq:helmholtz}) dependent
upon the remaining choice of boundary conditions for $\phi$ along
the swept out surface ${\bf \Gamma}$. For any such choice, we obtain
the gradient with respect to ${\bf W}$ as
\begin{equation}
\nabla_{{\bf W}}E={\bf V}^{\perp}={\bf V}-\nabla\phi\label{eq:gradient}
\end{equation}
to obtain a perturbation that maintains the solenoidal constraint
over ${\bf \Omega}$ while also maintaining the initial and final
densities. Plugging this into (\ref{eq:variation}) eliminates the
dependency upon the Lagrange multiplier $\lambda$ as well as the
integrals along the temporal faces ${\bf \Gamma}_{0}$ and ${\bf \Gamma}_{1}$
(see Appendix \ref{app:variation}), yielding 
\begin{equation}
\delta E\Bigr|_{_{\delta{\bf W}={\bf V}^{\perp}}}=\int_{{\bf \Omega}}\|{\bf V}-\nabla\phi\|^{2}\;d{\bf X}+\int_{{\bf \Gamma}}\left(\phi+\alpha\right)\,\left({\bf V}-\nabla\phi\right)\cdot{\bf {\bf N}}+{\bf W}\cdot\left({\bf V}+\nabla_{S}\alpha\right)\delta{\bf \Gamma}\cdot{\bf N}\,dS\label{eq:gradient_variation}
\end{equation}

\subsubsection{Constrained gradient with fixed support (Neumann conditions)}

While the Helmholtz decomposition (\ref{eq:helmholtz}) is not unique
if the only criterion is to split the vector field into purely irrotational
and solenoidal contributions, uniqueness can be obtained in a special
case by seeking contributions that are also orthogonal. We observe
that

\begin{align*}
\int_{{\bf \Omega}}{\bf V}^{\|}\cdot{\bf V}^{\perp}d{\bf X} & =\int_{{\bf \Omega}}\underbrace{\nabla\phi}_{{\bf V}^{\parallel}}\cdot\underbrace{\left({\bf V}-\nabla\phi\right)}_{{\bf V}^{\perp}}\,d{\bf X}=\underbrace{-\int_{{\bf \Omega}}\phi\left(\nabla\cdot{\bf V}-\Delta\phi\right)\,d{\bf X}}_{0\mbox{ for any solution }\Delta\phi=\nabla\cdot{\bf V}}+\underbrace{\int_{\partial{\bf \Omega}}\phi\left({\bf V}-\nabla\phi\right)\cdot{\bf N}\,d{\bf X}}_{0\mbox{ for Neumann }\nabla\phi\cdot{\bf N}={\bf V}\cdot{\bf N}}
\end{align*}
and notice that the final region integral over ${\bf \Omega}$ disappears
for any solution of the Poisson equation (\ref{eq:poisson}), regardless
of the choice of boundary condition. However, to obtain orthogonality
the additional boundary integral term above must also disappear. This
happens by imposing Neumann conditions
\begin{equation}
\nabla\phi\cdot{\bf N}={\bf V}\cdot{\bf N},\;{\bf X}\in{\bf \Gamma}\label{eq:neumann}
\end{equation}
(the same type of Neumann conditions already imposed along the temporal
boundaries ${\bf \Gamma}_{0}$ and ${\bf \Gamma}_{1}$) and yields
a gradient ${\bf V}^{\perp}$ with vanishing flux ${\bf V}^{\perp}\cdot{\bf N}=0$
along all boundaries of the support.

However, it also constrains the normal perturbation $\delta{\bf \Gamma}\cdot{\bf N}=0$
of the boundary itself, which is coupled to the flux perturbation
${\bf \delta}{\bf W}\cdot{\bf N}$ as follows (see Appendix \ref{app:coupled})
\begin{equation}
\delta{\bf W}\cdot{\bf N}=\sum_{i=1}^{n}\frac{\partial}{\partial s_{i}}\Bigl(({\bf W}\cdot{\bf T}_{i})(\delta{\bf \Gamma}\cdot{\bf N})\Bigr)\label{eq:coupling}
\end{equation}
where $s_{1},\ldots,s_{n}$ denote isothermal coordinates for ${\bf \Gamma}$
along the principal directions ${\bf T}_{1},\ldots,{\bf T}_{n}$ (unit
tangent representations). In particular, the vanishing flux condition
${\bf \delta}{\bf W}\cdot{\bf N}=0$ causes the coupled variation
of the support boundary to vanish in the unit normal direction as
well, $\delta{\bf \Gamma}\cdot{\bf N}=0$. Both of these effects cause
the boundary integral terms in (\ref{eq:gradient_variation}) to drop
away, independent of the remaining Lagrange multiplier $\alpha$,
thereby yielding
\begin{align*}
\left.\delta E\right|_{_{\delta{\bf W}={\bf V}^{\perp}}}= & \int_{{\bf \Omega}}\|{\bf V}^{\perp}\|^{2}\,d{\bf X}
\end{align*}
\emph{The key point to make here is that imposing a zero flux perturbation
along the original support boundary during the optimization process,
automatically constrains the evolution (and therefore the final optimizer)
to keep the original support. If we do not wish to constrain the solution
in this way, then we cannot impose vanishing flux conditions, but
as shown next, we must replace the Neumann conditions with strategically
chosen Dirichlet conditions instead. This will allow the generation
of inward or outward flux along the original boundary, which provides
information on how the boundary itself should evolve according to
the coupling in (\ref{eq:coupling}).}

\subsubsection{Unconstrained gradient with evolving support (Dirichlet conditions)}

To obtain the unconstrained gradient, in which the support boundary
is also allowed to evolve as part of the optimization process, we
must choose the remaining Lagrange multiplier as well as the boundary
conditions for $\phi$ to eliminate the boundary integral terms in
(\ref{eq:gradient_variation}) even when the coupled flux perturbation
$\delta{\bf W}\cdot{\bf N}$ and boundary perturbation ${\bf \delta{\bf \Gamma}\cdot{\bf N}}$
are both non-zero. We begin by solving the geometric transport PDE
\begin{equation}
{\bf W}\cdot\left({\bf V}+\nabla_{S}\alpha\right)=0\label{eq:geometric-transport}
\end{equation}
for the Lagrange multiplier $\alpha$ along the hypersurface ${\bf \Gamma}$
to remove the sensitivity of (\ref{eq:gradient_variation}) with respect
to $\delta{\bf \Gamma}$. We may reformulate this into a more standard
volumetric linear transport PDE over the full spatio-temporal domain
${\bf \Omega}$ by solving for a differentiable extension of $\alpha$
where we may express the intrinsic gradient $\nabla_{S}\alpha$ along
the boundary as the orthogonal projection of the gradient $\nabla\alpha$
of the volumetric extension as follows

\[
\nabla_{S}\alpha=\left({\cal I}-{\bf N}{\bf N}^{T}\right)\nabla\alpha
\]
Substituting this into the previous equation yields 
\begin{align*}
{\bf W}\cdot{\bf V}+{\bf W}^{\cdot T}\left({\cal I}-{\bf N}{\bf N}^{T}\right)\nabla\alpha & =0\\
{\bf W}\cdot{\bf V}+{\bf W}^{\cdot}\cdot\nabla\alpha-\underbrace{\left({\bf W}^{\cdot}\cdot{\bf N}\right)}_{0}\left(\nabla\alpha\cdot{\bf N}\right) & =0
\end{align*}
demonstrating that the solution for $\alpha$ along the boundary ${\bf \Gamma}$
does not depend upon its resulting extension ${\bf \nabla\alpha\cdot{\bf N}}$
along the normal due to the vanishing flux property ${\bf W}\cdot{\bf N}=0$.
As such, we obtain the following non-homogeneous linear transport
equation for $\alpha$
\begin{equation}
{\bf W}\cdot\nabla\alpha=-{\bf W}\cdot{\bf V}\label{eq:linear-transport}
\end{equation}
which is easily solved volumetrically and whose solution along the
boundary yields the desired solution for (\ref{eq:geometric-transport}).

After solving (\ref{eq:linear-transport}) for $\alpha$, we then
set $\phi=-\alpha$ along the boundary to eliminate the sensitivity
of (\ref{eq:gradient_variation}) with respect to the normal derivative
$\nabla\phi\cdot{\bf N}$ along the boundary, thereby transforming
the Neumann boundary condition into a Dirichlet condition instead.
\begin{equation}
\phi=-\alpha,\;{\bf X}\in{\bf \Gamma}\label{eq:dirichlet}
\end{equation}
Solving the complete system (\ref{eq:poisson}), (\ref{eq:temporal_bc}),
and (\ref{eq:dirichlet}) yields an unconstrained gradient descent
perturbation

\[
\delta{\bf W}=-{\bf V}^{\perp}=\nabla\phi-{\bf V}
\]
with a non-zero flux perturbation $\delta{\bf W}\cdot{\bf N}\ne0$.

\subsubsection{Initialization and gradient descent}

We now show how any initial choice of swept out hypersurface ${\bf \Gamma}_{\mbox{init}}$
and solenoidal vector field ${\bf W}_{\mbox{init}}$ may be deformed
using gradient descent within the class of smooth boundary surfaces
and solenoidal vector fields in order to solve the compactly supported
optimal transport problem directly in $\mathbb{R}^{n+1}$. 

\paragraph{Computing an initial solenoidal field ${\bf W}$}

If we combine the solenoidal constraint and boundary flux conditions
for ${\bf W}$ with the additional constraint that the initial field
${\bf W}_{\mbox{init}}$ be conservative as well, then we may plug
${\bf W=\nabla{\bf \Phi}}$ into (\ref{eq:solenoidal}) and (\ref{eq:flux}),
for some scalar spatiotemporal function ${\bf \Phi}$, to obtain the
Laplace equation with Neumann boundary conditions (non-homogeneous
along the two temporal faces).
\begin{align}
\Delta{\bf \Phi} & =0,\qquad\quad\text{{\bf X}}\in{\bf \Omega}\label{eq:harmonic}\\
\nabla{\bf \Phi}\cdot{\bf N} & =\begin{cases}
\;\;0\;\;, & {\bf X}\in{\bf \Gamma}\\
-\rho_{0}, & {\bf X}\in{\bf \Gamma}_{0}\\
+\rho_{1}, & {\bf X}\in{\bf \Gamma}_{1}
\end{cases}\nonumber 
\end{align}
A solution will exist so long as $0=\int_{d{\bf \Omega}}\nabla{\bf \Phi}\cdot{\bf N}\,d{\bf S}=\int_{d{\bf \Omega}}{\bf W}\cdot{\bf N}\,d{\bf S}$,
which in this case is equivalent to our balanced assumption $\int_{\Omega_{0}}\rho_{0}\,dx=\int_{1}\rho_{1}\,dx$
for the initial and target densities. The solution will be unique
up to an additive constant for ${\bf \Phi}$, which will then disappear
after applying the gradient to obtain the following initial vector
field.
\begin{equation}
{\bf W}_{\mbox{init}}=\nabla{\bf \Phi}\label{eq:conservative}
\end{equation}

\paragraph{Gradient step for the solenoidal field ${\bf W}$}

A descent step on ${\bf W}$ can be taken by computing the gradient
(\ref{eq:gradient}) either for the fixed support case by solving
the Poisson equation (\ref{eq:poisson}) with the Neumann boundary
conditions (\ref{eq:neumann}) or, for the case of co-evolving support,
using the Dirichlet boundary conditions (\ref{eq:dirichlet}) obtained
through (\ref{eq:linear-transport}). We may then take a gradient
descent step
\begin{align}
{\bf W}_{k+1} & \to{\bf W}_{k}-\gamma_{k}{\bf V}_{k}^{\perp}\label{eq:descent_step}
\end{align}
where ${\bf V}_{k}^{\perp}={\bf V}_{k}-\nabla\phi$. Initially, after
the initialization strategy outlined above, use of the Neumann condition
to first optimize with respect to ${\bf W}$ over the initially chosen
support, is recommend prior to switching to the joint evolution Dirichlet
strategy. During this initial optimization step with fixed support,
we may use Newton's method to determine the optimal step factor $\gamma_{k}$
by solving (see Appendix \ref{sec:newton})

\begin{align*}
0 & =\frac{dE_{k+1}}{d\gamma_{k}}=\int_{{\bf \Omega}}{\bf V}_{k+1}(\gamma_{k})\cdot\left(\nabla\phi_{k}-{\bf V}_{k}\right)\,d{\bf X}
\end{align*}
for each gradient step.

\paragraph{Gradient step for the surface ${\bf \Gamma}$}

When we apply the Dirichlet update to ${\bf W}$, the original flux-free
solenoidal field will develop non-vanishing flux ${\bf W}\cdot{\bf N}$
along the original boundary. If we now change the Lagrange multiplier
\[
\alpha=-\phi
\]
to eliminate the sensitivity of (\ref{eq:gradient_variation}) with
respect to the flux evolution $\delta{\bf W}\cdot{\bf N}$, we obtain
the combined sensitivity with respect to both the evolving solenoidal
filed ${\bf W}$ and its support as follows.
\begin{equation}
\delta E\Bigr|_{_{\delta{\bf W}={\bf V}^{\perp}}}=\int_{{\bf \Omega}}\|{\bf V}-\nabla\phi\|^{2}\;d{\bf X}+\int_{{\bf \Gamma}}{\bf W}\cdot\left({\bf V}-\nabla\phi\right)\delta{\bf \Gamma}\cdot{\bf N}\,dS\label{eq:gradient_boundary_variation}
\end{equation}
The contribution for the first integral term has already been established
in our solution of the Dirichlet problem for $\phi$. To maximize
the contribution of the remaining surface integral term, we apply
the following gradient perturbation to ${\bf \Gamma}$.
\begin{equation}
\delta{\bf \Gamma}={\bf W}\cdot\left({\bf V}-\nabla\phi\right){\bf N}\label{eq:gradient_boundary}
\end{equation}

\section{Preliminary experimental results}

We conclude with two experimental results which illustrate the benefits
of this variational approach, stemming in particular from its separate
yet coupled optimization of the compact spatiotemporal support and
the density within. While the mathematical formulation of the approach
has been fully outlined here, the numerical implementation strategies
are still under investigation and, as such, the following results
are intended to be preliminary indications of what we may expect after
more sophisticated numerical strategies have been further explored
and developed.

\subsection{Interpolation between two different non-convex supports}

\begin{figure}[ht]
\begin{tabular}{c@{}c@{}c@{}c}
\vspace{-0.5mm}
real image & real image & real image & real image \\
\epsfig{figure=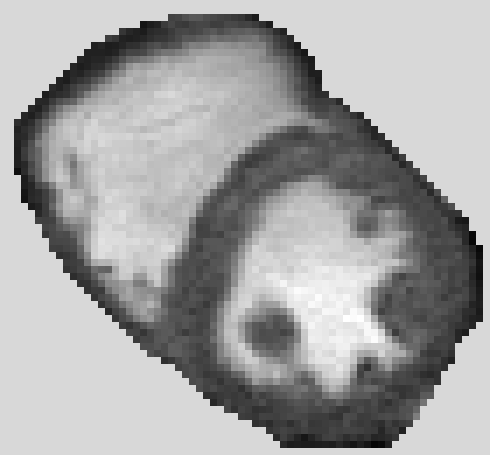,width=1.6in} &
\epsfig{figure=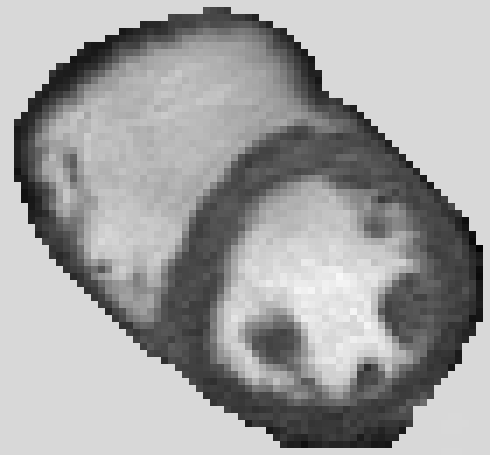,width=1.6in} &
\epsfig{figure=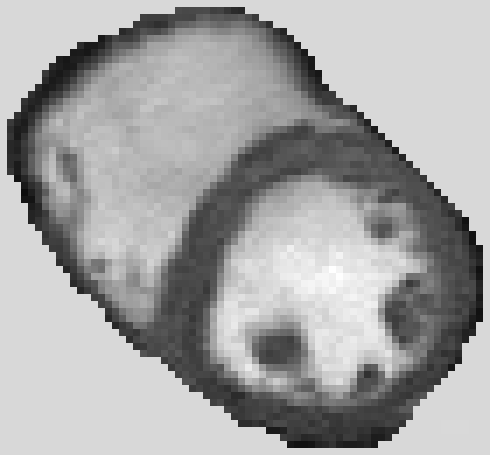,width=1.6in} &
\epsfig{figure=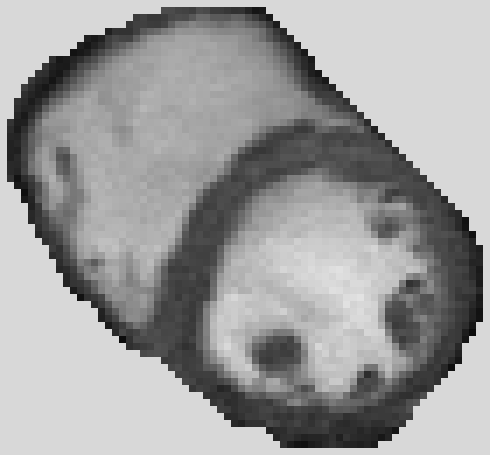,width=1.6in} \\
$t=0$ & $t=\frac{1}{3}$ & $t=\frac{2}{3}$ & $t=1$ \\
\epsfig{figure=density0,width=1.6in} &
\epsfig{figure=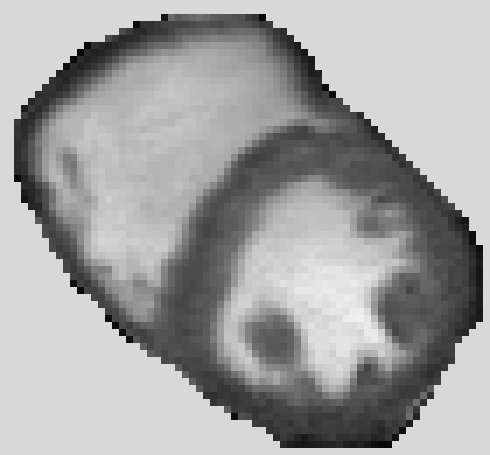,width=1.6in} &
\epsfig{figure=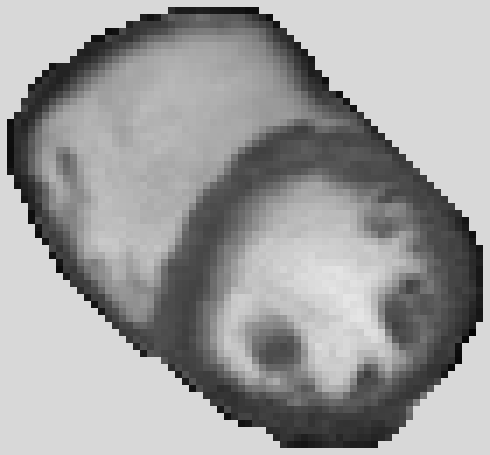,width=1.6in} &
\epsfig{figure=density3,width=1.6in}
\vspace{-1.5mm}
\\
real image & transport image & transport image & real image
\end{tabular} \vspace{-5mm}
\caption{Density evolution between two non-convex, differing compact supports.
Left and right input images ($t$=0 and $t$=1) show the
starting and ending densities and supports (gray background does not
represent any density value and has no effect on the computations),
while the middle images (bottom) represent transported densities and supports
computed at equally spaced intermediate times (with the actual measured cardiac
images at corresponding times shown along the top row for comparison).
\label{fig:cariac-density}
}
\end{figure}

\begin{wrapfigure}[24]{o}{0.31\columnwidth}
\vspace{-10mm}
\epsfig{figure=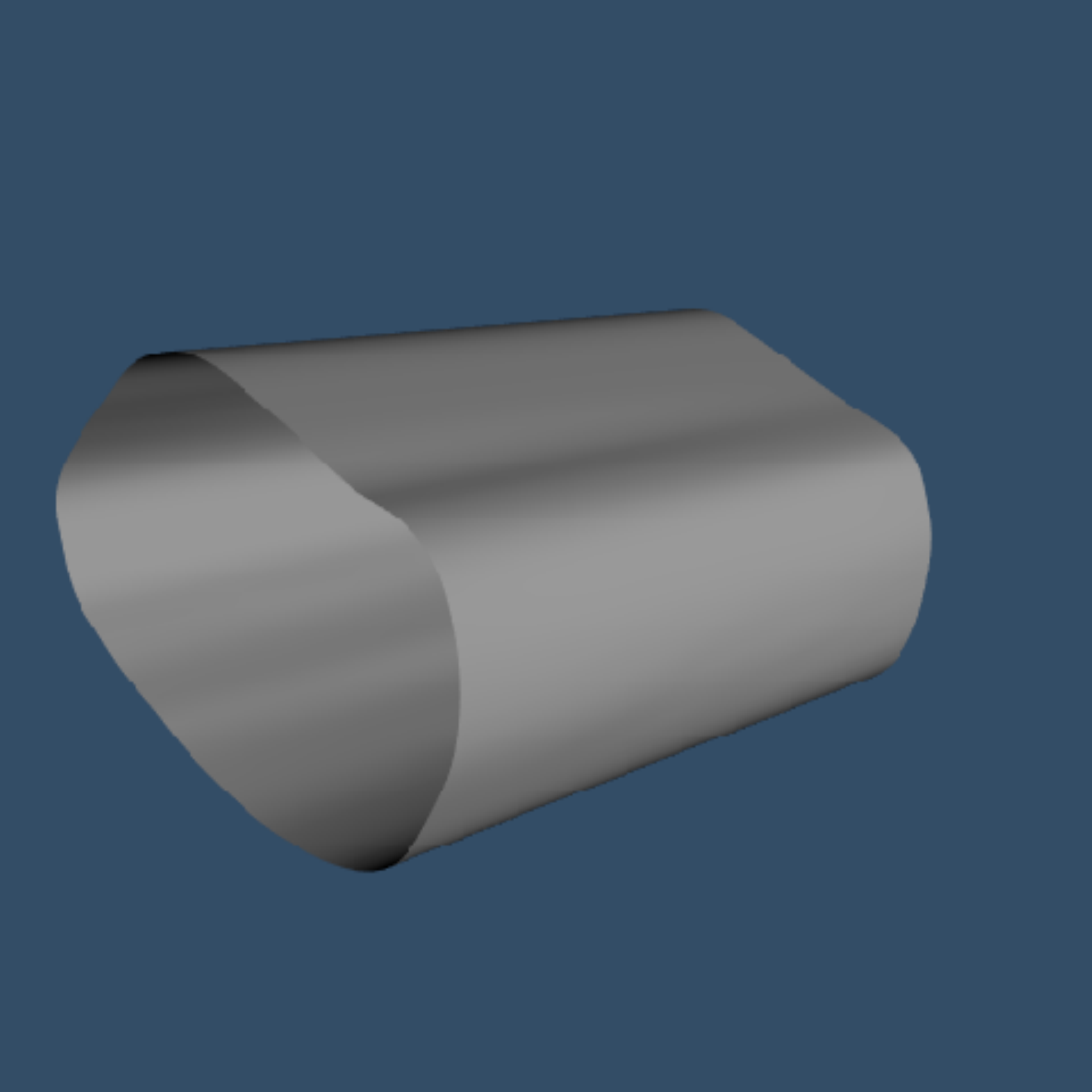,width=2in}\vspace{-10mm}
\epsfig{figure=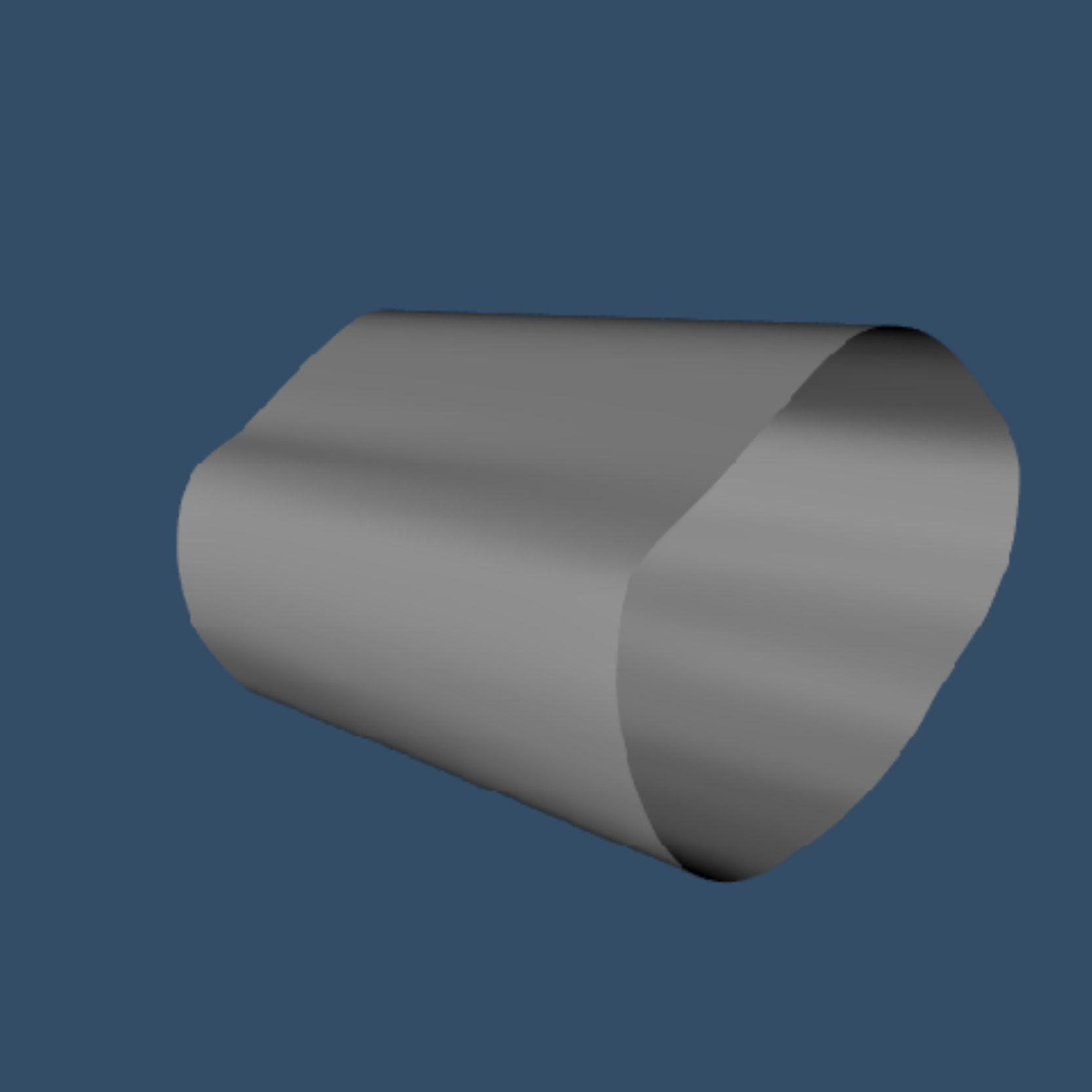,width=2in}\vspace{-3.5mm}
\caption{Two rendered viewpoints of the hypersurface ${\bf\Gamma}$ which
models/constrains the dynamic support from $t$=0 to $t$=1.
\label{fig:cariac-support}}
\end{wrapfigure}

In this first example we tackle the problem of interpolating between
cardiac images captured at two different moments within the heart beat
cycle shown on the left and the right in Figure \ref{fig:cariac-density}
Notice that both cell boundaries represent non-convex shapes with
several small concavities. We also see structures of
interest inside the ventricle (papillary muscle cross sections)
which not only move and
deform along with the rest of the image, but which also change in their
topological appearance.
At a coarse scale, the boundary shapes appear to be similar, but detailed
inspection reveals that finer scale protrusions and intrusions
around the boundary differ, especially the concavity along the upper right
edge of the left image which disappears in the right image.
Nevertheless, both boundaries exhibit a matching
simple topology, which we would like to preserve when morphing one
into the other. Numerically, this can be challenging without an explicit
model of the support boundary, making it difficult to guarantee that
small scale protrusions don't break off during transport to yield
transitional topological changes. 

Using these two cardiac images as the starting and ending densities
at time $t=0$ and $t=1$ respectively, we solve the compact optimal
transport problem with the variational approach outlined in this paper,
using a 3D space-time grid for the solenoidal vector field ${\bf W}$
with 64 temporal slices, each of size 128x128 (same resolution as
the two input images). A matching level set grid $\Psi$ is used to
represent the spatiotemporal support as the set $\Psi<0$. We see
the density (temporal component of ${\bf W}$) in Figure \ref{fig:cariac-density}
at equally space intervals along the computed transport.
We can see this more explicitly by visualizing the entire swept-out
hypersurface ${\bf \Gamma}$ (the portion of the spatiotemporal support
boundary strictly between 0 and 1) in Figure \ref{fig:cariac-support}.
This is very easily rendered as the zero level set of $\Psi$ and
clearly reveals a smooth homotopy connecting the two end curves, one
exposed within the top rendering and the other within the bottom).

\subsection{Optimal transport with spatiotemporal support constraints}

This next experiment illustrates through and intuitive toy-example
an important and useful extension to the class of optimal transport
problems which, to the best of our knowledge, is not accommodated
by prior art, but is easily handled within the methodology outlined
here. Namely, how might one compute the optimal transport between
densities subject to constraints imposed on where (and possibly even
when) mass is allowed or not allowed to move along the way? 

\subsubsection{Known initial and final densities with intermediate support constraints}

We can easily motivate the utility of such constraints if we go back
to the classic problem posed by Gaspard Monge. In formulating the
problem to optimally move a pile of dirt from one place to another,
no constraints were imposed on the path taken by each portion of moved
dirt. While the unconstrained optimal solution may yield a realistic
and realizable transport strategy in several practical circumstances,
this cannot always be guaranteed. For example, suppose the task is
to move a pile of dirt across a limited number of bridges to the other
side of a river. The unconstrained solution could easily yield an
impractical transport strategy which involves crossing open portions
of the river. A related larger scale problem might involve planning
the motion of massive numbers of land troops distributed over a set
of territories to another set of territories taking into account geographical
barriers such as mountains and bodies of water as well as political
barriers which would render certain intermediate territories out of
bounds.

Even when the topology of the desired transport is known, optimization
subject to geometric constraints can still be nontrivial. For example,
transporting a single pile of dirt across a single bridge which is
narrower than the base of the pile is already an interesting problem.
The unconstrained optimal transport will likely want to move some
portion of dirt outside the confines of the bridge. Should all of
that excess dirt simply be re-routed and accumulated along the closest
edge of the bridge, or should some of it be moved more centrally inside
the bridge, which makes the trajectory deviate even further from the
unconstrained optimum but attenuates an otherwise massive spike in
density along the edge? Clearly there is a trade off that is not easily
intuited directly from the unconstrained solution.

These types of constraints are easily imposed using the variational
framework presented in this paper due to the explicit representation
and separate-yet-coupled evolution of the support boundary. So long
as we choose an initial spatiotemporal support that satisfies the
provided set of spatiotemporal constraints, the calculated gradient
descent evolution of the resulting spatiotemporal hypersurface can
simply be set to zero locally wherever its application would otherwise
violate the constraints. 

\subsubsection{Unknown final density but with known support}

\begin{figure}[t]
\begin{tabular}{c@{}c@{}c@{}c}
\epsfig{figure=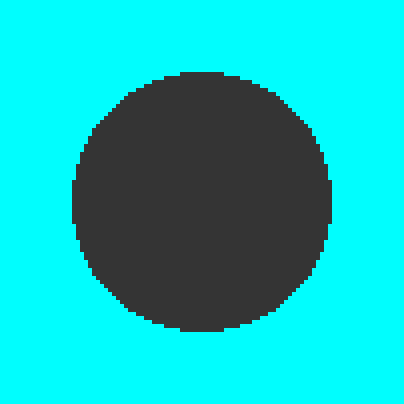,width=1.6in} \hspace{-8mm} &
\epsfig{figure=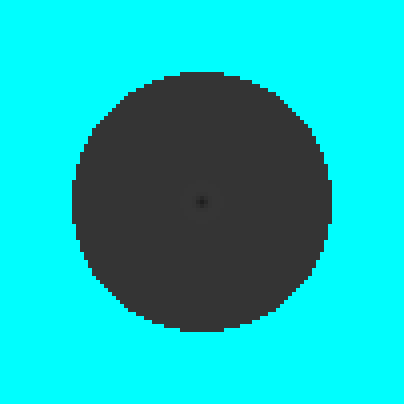,width=1.6in} \hspace{-8mm} &
\epsfig{figure=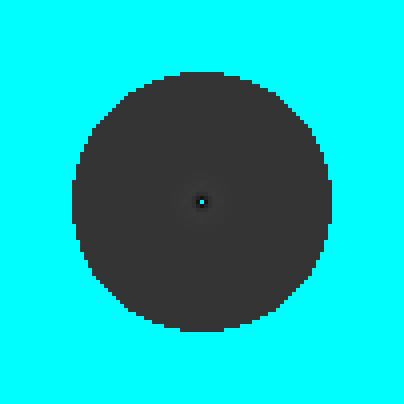,width=1.6in} \hspace{-8mm} &
\epsfig{figure=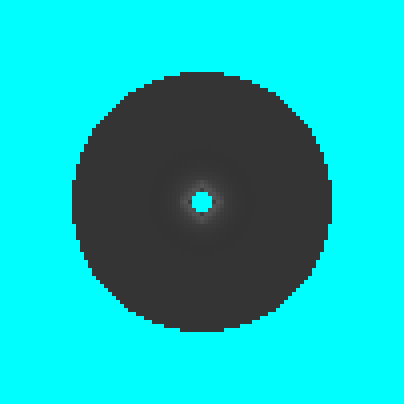,width=1.6in} \vspace{-8mm} \\
\epsfig{figure=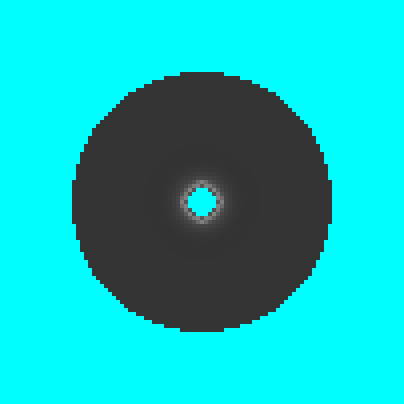,width=1.6in} \hspace{-8mm} &
\epsfig{figure=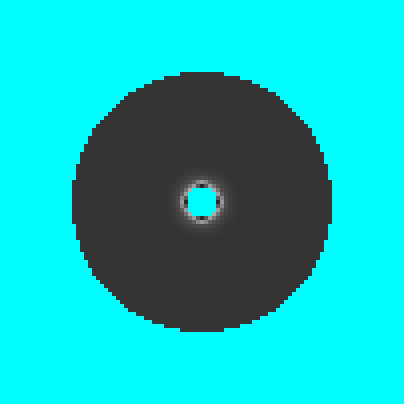,width=1.6in} \hspace{-8mm} &
\epsfig{figure=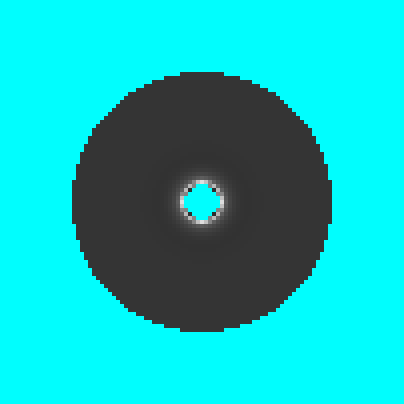,width=1.6in} \hspace{-8mm} &
\epsfig{figure=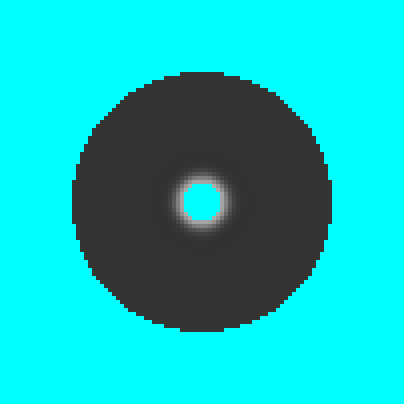,width=1.6in}
\end{tabular} \vspace{-4mm}
\caption{Example of optimal transport with spatiotemporal support constraints.
Mass is evacuated from the center of a disk with uniform initial density
of 3.0 (appears as dark-gray/almost-black in the top-left frame),
while constrained to remain inside the disk. The optimal evacuation
strategy is computed by solving the optimal transport problem between
the initial uniform density on the disk back to itself (normally the
trivial transport of zero) but imposing a hole in the spatiotemporal
support at $t$=$\frac{1}{2}$. Mass therefore evacuates the hole
from $t$=0 to $t$=$\frac{1}{2}$ (shown above) and refills the hole
from $t$=$\frac{1}{2}$ to $t$=1 (reverse of above).\label{fig:evacuate}}
\end{figure}

Another extension of problems that are easily accommodated by this
coupled variational approach include scenarios where the support of
the target density is given but the target density itself is unknown
(and therefore part of the optimization). Such a problem can easily
be transformed into the problem of a known final density with intermediate
constraints on the support by treating the desired final density as
the halfway point ($t=\frac{1}{2}$) in transporting the initial density
at $t=0$ back to itself again at $t=1$. In this way, the desired
final support becomes a constraint on the intermediate support instead.
Optimization of this reconfigured problem will yield both a forward
copy (from $t=0$ to $t=\frac{1}{2}$) and a backward copy (from $t=\frac{1}{2}$
to $t=1$) of the optimal transport for the original problem as well
as the optimal target density itself at $t=\frac{1}{2}$. As such,
from an implementation standpoint, this class of problems can be handled
the same way as the class of problems just described above.

A practical application for this form of constrained optimal transport
would be to compute the most efficient evacuation strategy to clearing
mass out of a given sub-region while keeping it within a larger encompassing
region that already contains pre-existing mass. In this case we know
the initial density and support, and we know the final support, simply
the initial support minus the sub-region to be evacuated, but we don't
know or otherwise want to constrain the resulting new density within
the now reduced support. We therefore seek the least costly way (according
to the Wasserstein metric) to redistribute the mass originally contained
within the subregion into its surrounding, already occupied neighborhood.
Such a problem may arise, for example, when seeking to clear extensive
zones of all materials and/or personnel while keeping them within
the confines of larger zones whose occupants are free to be internally
relocated if needed.

We illustrate precisely this scenario with an intuitive toy-problem
in Figure \ref{fig:evacuate}. We start with uniformly distributed
mass at $t=0$ with density 3.0 across a disk representing the global
confines. In turn, we define the target distribution at $t=1$ to
be the same, but impose a hole in the intermediate support at $t=\frac{1}{2}$
within the center of the disk. Solving the constrained optimal transport
problem between these matching uniform distributions causes the initially
filled hole to be evacuated from $t=0$ to $t=\frac{1}{2}$, as illustrated
from left-to-right and top-to-bottom in Figure \ref{fig:evacuate},
and then to be refilled from $t=\frac{1}{2}$ to $t=1$, as also illustrated
in Figure \ref{fig:evacuate} when read in reverse. The optimal redistribution
of mass is attained at $t=\frac{1}{2}$ and is displayed in grayscale
at the end of the figure.

\begin{wrapfigure}[30]{o}{4.3in}
\vspace{-4mm}
\begin{tabular}{|c|c|}
\hline 
\epsfig{figure=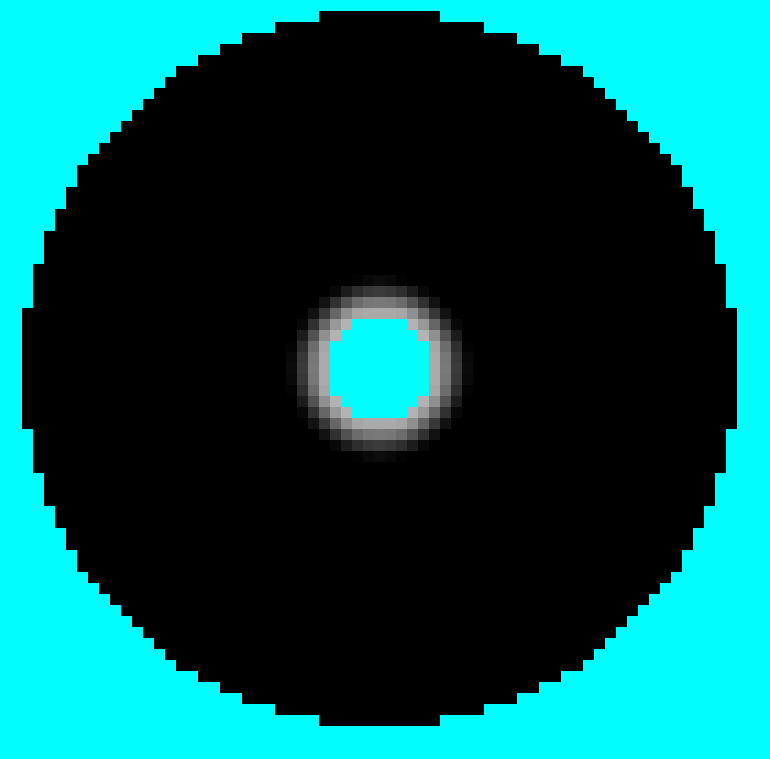,width=1.8in} &
\epsfig{figure=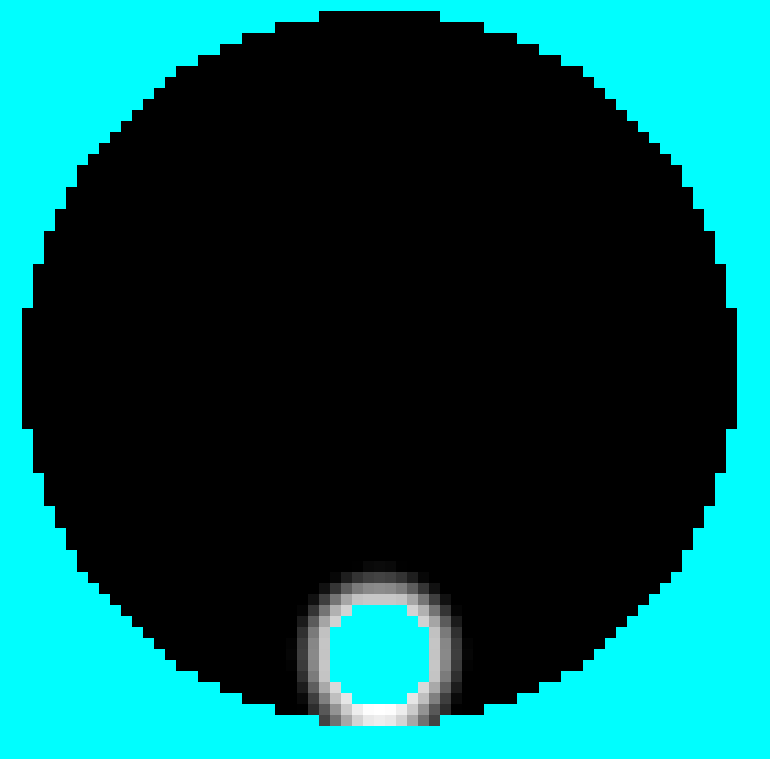,width=1.8in} \vspace{-1mm} \\
Evacuated mass from disk center &
Evacuated mass near disk edge \\
\epsfig{figure=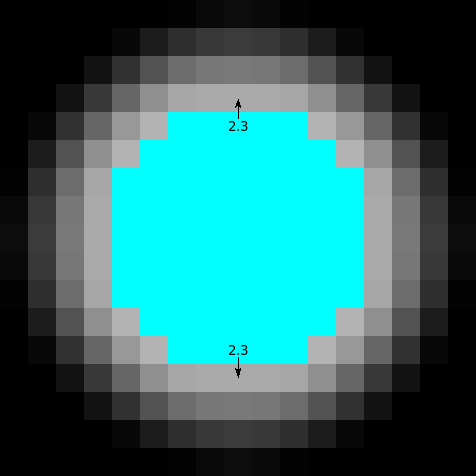,width=1.8in} &
\epsfig{figure=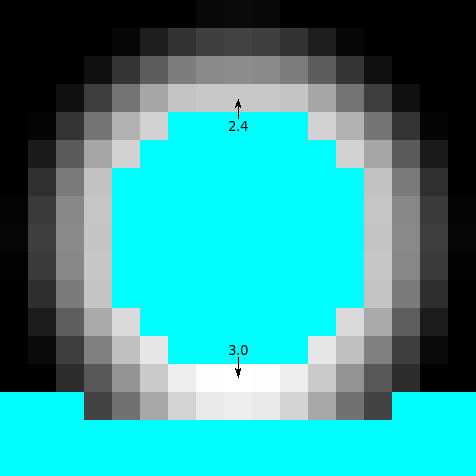,width=1.8in} \vspace{-1mm} \\
Close-up: symmetric increase &
Close-up: nonsymmetric increase \\
\hline 
\end{tabular} \vspace{-4mm}
\caption{Optimized mass evacuation from two different locations.
Both cases begin with uniform mass density (3.0) over the disk
then evacuate a hole. Accumulated density results
wherever mass gets relocated (see previous figure). Final net accumulation
is shown here as the rise in density.
When the hole is perfectly centered (left), evacuated
mass gets redistributed symmetrically with a peak density rise of 2.3 (77\%).
Non-symmetric restribution produces a peak rise of 3.0 when the
hole is created near the boundary (right).
\label{fig:accumulate}}
\end{wrapfigure}

Even in this simplest illustrative example with constant density and
maximal symmetry, it is by no means intuitively obvious how far away
mass should be displaced from the hole boundary compared to how much
it should be allowed to accumulate along the boundary. In fact, the
density would become infinite if the evacuated mass were to remain
strictly along the boundary. To get a better sense of where displaced\footnote{displaced mass includes not only evacuated mass from the hole, but
also some mass originally outside the hole which gets transferred
further away to make room for the evacuated mass} mass accumulates upon preexisting mass, we show the net density change
in Figure \ref{fig:accumulate} (left side) which attains a maximum
rise of 2.3 all around the boundary of the evacuated hole and gradually
rolls off further outward.

We can make this toy problem even more interesting by evacuating mass
near the boundary of the disk rather than its center. Looking at the
resulting net density increase in Figure \ref{fig:accumulate} (right
side), we can make several observations. First, as expected, the redistribution
is no longer symmetric. The symmetry is broken two different ways.
First, the rise in density is much higher (3.0, a full 100\% jump)
along the bottom border of the hole compared to the top border. This
is unsurprising since there is not as much room to move away from
the hole, and so the same amount of evacuated mass distributed over
a thinner local neighborhood must necessarily result in a larger accumulated
density.

However, we also notice that the density jump near the top
of the hole (2.4), while lower than the bottom, still exceeds the
symmetric jump (2.3) when the hole was centered inside the disk. This
differential grows as we travel along the border of the hole toward
the bottom. This means that some of the mass within the lower half
of the hole, which was all evacuated downward from the centered hole,
actually got evacuated upward from the hole near the boundary of the
confining domain. As such, at least two interesting interplays are
relevant in this optimization. First, as in the symmetric case, how
far away should mass be displaced away from borders versus accumulated
along borders? Second, when displacement distances are limited, what
is the right trade-off between higher accumulation at nearby borders
versus more costly redirection toward farther borders of the region
to be evacuated? These considerations are both naturally handled by
this coupled variational framework.

\bibliographystyle{ieeetr}
\bibliography{refs}

\begin{thebibliography}{10}

\bibitem{GM}
G.~Monge, ``M{\'e}moire sur la th{\'e}orie des d{\'e}blais et des remblais,''
  {\em Histoire de l'Acad{\'e}mie Royale des Sciences de Paris}, 1781.

\bibitem{LK}
L.~V. Kantorovich, ``On a problem of monge,'' {\em CR (Doklady) Acad. Sci. URSS
  (NS)}, vol.~3, pp.~225--226, 1948.

\bibitem{GangboMcCann}
W.~Gangbo and R.~McCann, ``The geometry of optimal transportation,'' {\em Acta
  Mathematica}, vol.~177, no.~2, pp.~113--161, 1996.

\bibitem{villani1}
C.~Villani, {\em Topics in Optimal Transportation}.
\newblock American Mathematical Soc., 2003.

\bibitem{villani2}
C.~Villani, {\em Optimal Transport: Old and New}, vol.~338.
\newblock Springer Science \& Business Media, 2008.

\bibitem{Evans1999}
L.~C. Evans and W.~Gangbo, ``{Differential equations methods for the
  Monge-Kantorovich mass transfer problem},'' {\em Memoirs of the American
  Mathematical Society}, vol.~137, no.~653, pp.~0--0, 1999.

\bibitem{Rachev1998}
S.~T. Rachev and L.~R{\"{u}}schendorf, {\em {Mass Transportation Problems:
  Volume I: Theory}}.
\newblock Probability and its Applications, Berlin: Springer, 1998.

\bibitem{Arjovsky2017}
M.~Arjovsky, S.~Chintala, and L.~Bottou, ``{Wasserstein GAN},'' {\em
  arxiv.org}, vol.~1701.07875, 2017.

\bibitem{CarMaa14}
E.~A. Carlen and J.~Maas, ``An analog of the 2-{W}asserstein metric in
  non-commutative probability under which the {F}ermionic {F}okker--{P}lanck
  equation is gradient flow for the entropy,'' {\em Communications in
  Mathematical Physics}, vol.~331, no.~3, pp.~887--926, 2014.

\bibitem{Haker2004}
S.~Haker, L.~Zhu, A.~Tannenbaum, and S.~Angenent, ``{Optimal mass transport for
  registration and warping},'' {\em International Journal of Computer Vision},
  vol.~60, no.~3, pp.~225--240, 2004.

\bibitem{MitMie16}
M.~Mittnenzweig and A.~Mielke, ``An entropic gradient structure for lindblad
  equations and generic for quantum systems coupled to macroscopic models,''
  {\em arXiv preprint arXiv:1609.05765}, 2016.

\bibitem{rr}
S.~T. Rachev and L.~R{\"u}schendorf, {\em Mass Transportation Problems:
  {V}olumes {I} and {II}}.
\newblock Springer Science \& Business Media, 1998.

\bibitem{Statement2020}
J.~C. Mathews, S.~Nadeem, M.~Pouryahya, Z.~Belkhatir, J.~O. Deasy, A.~J.
  Levine, and A.~R. Tannenbaum, ``{Functional network analysis reveals an
  immune tolerance mechanism in cancer},'' {\em Proceedings of the National
  Academy of Sciences}, vol.~117, pp.~16339--16345, jul 2020.

\bibitem{BB}
J.-D. Benamou and Y.~Brenier, ``A computational fluid mechanics solution to the
  {M}onge-{K}antorovich mass transfer problem,'' {\em Numerische {M}athematik},
  vol.~84, no.~3, pp.~375--393, 2000.

\bibitem{Otto}
F.~Otto, ``The geometry of dissipative evolution equations: the porous medium
  equation,'' {\em Communications in Partial Differential Equations}, vol.~26,
  no.~1-2, pp.~101--174, 2001.

\bibitem{omtEx3}
L.~Zhu, Y.~Yang, S.~Haker, and A.~Tannenbaum, ``An image morphing technique
  based on optimal mass preserving mapping,'' {\em IEEE Transactions on Image
  Processing}, vol.~16, no.~6, pp.~1481--1495, 2007.

\bibitem{haker}
S.~Haker, A.~Tannenbaum, and R.~Kikinis, ``Mass preserving mappings and image
  registration,'' {\em MICCAI}, pp.~120--127, 2001.

\bibitem{osher}
S.~Osher and R.~Fedkiw, {\em Level Set Methods and Dynamic Implicit Surfaces}.
\newblock Springer, 2003.

\bibitem{sethian}
J.~Sethian, {\em Level Set Methods Fast Marching Methods}.
\newblock Cambridge University Press, 1999.

\bibitem{sapiro}
G.~Sapiro, {\em Geometric Partial Equations and Image Analysis}.
\newblock Cambridge University Press, 2001.

\bibitem{susskind}
L.~Susskind and A.~Friedman, {\em Special Relativity and Classical Field
  Theory: The Theoretical Minimum}.
\newblock Hachette Book Group, 2017.

\end{thebibliography}

\section{Appendices}

\subsection{Detailed first variation and gradient calculation\label{app:variation}}

We first note that a variation of ${\bf W}$ yields a variation of
${\bf V}$ which is orthogonal to ${\bf W}$ itself. This may be demonstrated
directly using the expression (\ref{eq:V_W}) as follows.
\begin{align*}
\delta{\bf V} & =\left(\frac{\delta{\bf W}}{\left({\bf W}\cdot{\bf e}_{0}\right)}-\frac{{\bf W}}{\left({\bf W}\cdot{\bf e}_{0}\right)^{2}}\left(\delta{\bf W}\cdot{\bf e}_{0}\right)\right)-\left(\frac{\delta{\bf W}\cdot{\bf W}}{\left({\bf W}\cdot{\bf e}_{0}\right)^{2}}-\frac{{\bf W}\cdot{\bf W}}{\left({\bf W}\cdot{\bf e}_{0}\right)^{3}}\left(\delta{\bf W}\cdot{\bf e}_{0}\right)\right){\bf e}_{0}\\
{\bf W}\cdot\delta{\bf V} & =\underbrace{\left(\frac{\delta{\bf W}\cdot{\bf W}}{\left({\bf W}\cdot{\bf e}_{0}\right)}-\frac{{\bf W}\cdot{\bf W}}{\left({\bf W}\cdot{\bf e}_{0}\right)^{2}}\left(\delta{\bf W}\cdot{\bf e}_{0}\right)\right)-\left(\frac{\delta{\bf W}\cdot{\bf W}}{\left({\bf W}\cdot{\bf e}_{0}\right)}-\frac{{\bf W}\cdot{\bf W}}{\left({\bf W}\cdot{\bf e}_{0}\right)^{2}}\left(\delta{\bf W}\cdot{\bf e}_{0}\right)\right)}_{\mbox{terms cancel}}=0
\end{align*}
Using Lagrange multipliers to incorporate the mass preservation constraints
we write the energy as
\[
E({\bf W},{\bf \Gamma},\lambda,\alpha)=\int_{{\bf \Omega}}{\bf W}\cdot{\bf V}\,+\lambda\,\nabla\cdot{\bf W}\,d{\bf X}+\int_{{\bf \Gamma}}\alpha{\bf W}\cdot{\bf N}\,dS+\int_{{\bf \Gamma}_{0}}\alpha\left(\overbrace{{\bf W}\cdot{\bf N}}^{-{\bf W}\cdot{\bf {\bf e}_{0}}}+\rho_{0}\right)\,dS+\int_{{\bf \Gamma}_{1}}\alpha\left(\overbrace{{\bf W}\cdot{\bf N}}^{{\bf W}\cdot{\bf e}_{0}}-\rho_{1}\right)\,dS
\]
and compute its first variation as follows. 

\begin{align*}
\delta E= & \int_{{\bf \Omega}}\delta{\bf W}\cdot{\bf V}+\underbrace{{\bf W}\cdot\delta{\bf V}}_{0}+\lambda\nabla\cdot\delta{\bf W}+\delta\lambda\,\nabla\cdot{\bf W}\;d{\bf X}+\int_{{\bf \Gamma}}\left({\bf W}\cdot{\bf V}+\lambda\nabla\cdot{\bf W}\right)\delta{\bf \Gamma}\cdot{\bf N}\,dS\\
 & +\int_{{\bf \Gamma}}\delta\alpha\,{\bf W}\cdot{\bf N}+\alpha\,\delta{\bf W}\cdot{\bf N}\,dS+\int_{{\bf \Gamma}}\left(\nabla_{S}\alpha\cdot{\bf W}+\alpha\nabla\cdot{\bf W}\right)\delta{\bf \Gamma}\cdot{\bf N}\,dS\\
 & -\int_{{\bf \Gamma}_{0}}\delta\alpha\left({\bf W}\cdot{\bf e}_{0}-\rho_{0}\right)+\alpha\,\delta{\bf W}\cdot{\bf e}_{0}\,dS+\int_{{\bf \Gamma}_{1}}\delta\alpha\left({\bf W}\cdot{\bf e}_{0}-\rho_{1}\right)+\alpha\,\delta{\bf W}\cdot{\bf e}_{0}\,dS\\
= & \int_{{\bf \Omega}}\delta{\bf W}\cdot\left({\bf V}-\nabla\lambda\right)+\delta\lambda\,\nabla\cdot{\bf W}\;d{\bf X}+\int_{\partial{\bf \Omega}}\lambda\,\delta{\bf W}\cdot{\bf N}\,dS\\
 & +\int_{{\bf \Gamma}}\delta\alpha\,{\bf W}\cdot{\bf N}+\alpha\,\delta{\bf W}\cdot{\bf N}+\left({\bf W}\cdot\left({\bf V}+\nabla_{S}\alpha\right)+\left(\lambda+\alpha\right)\nabla\cdot{\bf W}\right)\delta{\bf \Gamma}\cdot{\bf N}\,dS\\
 & -\int_{{\bf \Gamma}_{0}}\delta\alpha\left({\bf W}\cdot{\bf e}_{0}-\rho_{0}\right)+\alpha\,\delta{\bf W}\cdot{\bf e}_{0}\,dS+\int_{{\bf \Gamma}_{1}}\delta\alpha\left({\bf W}\cdot{\bf e}_{0}-\rho_{1}\right)+\alpha\,\delta{\bf W}\cdot{\bf e}_{0}\,dS\\
\\
= & \int_{{\bf \Omega}}\delta{\bf W}\cdot\left({\bf V}-\nabla\lambda\right)+\delta\lambda\,\nabla\cdot{\bf W}\;d{\bf X}\\
 & +\int_{{\bf \Gamma}}\delta\alpha\,{\bf W}\cdot{\bf N}+\left(\lambda+\alpha\right)\,\delta{\bf W}\cdot{\bf N}+\left({\bf W}\cdot\left({\bf V}+\nabla_{S}\alpha\right)+\left(\lambda+\alpha\right)\nabla\cdot{\bf W}\right)\delta{\bf \Gamma}\cdot{\bf N}\,dS\\
 & -\int_{{\bf \Gamma}_{0}}\delta\alpha\left({\bf W}\cdot{\bf e}_{0}-\rho_{0}\right)+\left(\lambda+\alpha\right)\,\delta{\bf W}\cdot{\bf e}_{0}\,dS+\int_{{\bf \Gamma}_{1}}\delta\alpha\left({\bf W}\cdot{\bf e}_{0}-\rho_{1}\right)+\left(\lambda+\alpha\right)\,\delta{\bf W}\cdot{\bf e}_{0}\,dS
\end{align*}
Plugging in the mass conversation constraints eliminates the dependence
on $\delta\lambda$ and $\delta\alpha$ yielding the simpler expression
\begin{align*}
\delta E= & \int_{{\bf \Omega}}\delta{\bf W}\cdot\left({\bf V}-\nabla\lambda\right)\;d{\bf X}+\int_{{\bf \Gamma}}\left(\lambda+\alpha\right)\,\delta{\bf W}\cdot{\bf N}+{\bf W}\cdot\left({\bf V}+\nabla_{S}\alpha\right)\delta{\bf \Gamma}\cdot{\bf N}\,dS\\
 & -\int_{{\bf \Gamma}_{0}}\left(\lambda+\alpha\right)\,\delta{\bf W}\cdot{\bf e}_{0}\,dS+\int_{{\bf \Gamma}_{1}}\left(\lambda+\alpha\right)\,\delta{\bf W}\cdot{\bf e}_{0}\,dS
\end{align*}

Now apply the Helmholtz decomposition
\[
{\bf V}={\bf V}^{\parallel}+{\bf V}^{\perp}
\]
where

\begin{align*}
{\bf V}^{\parallel} & =\nabla\phi\\
{\bf V}^{\perp} & ={\bf V}-\nabla\phi\\
\Delta\phi & =\nabla\cdot{\bf V}\;\mbox{inside }{\bf \Omega}\quad\left(\mbox{\ensuremath{\therefore\nabla\cdot{\bf V}^{\perp}=0}}\right)\\
\nabla\phi\cdot\underbrace{{\bf N}}_{-{\bf e}_{0}} & ={\bf V}\cdot{\bf N}\;\mbox{along }{\bf \Gamma}_{0}\quad\left(\mbox{\ensuremath{\therefore{\bf V}^{\perp}\cdot{\bf e}_{0}=0}}\right)\\
\nabla\phi\cdot\underbrace{{\bf N}}_{+{\bf e}_{0}} & ={\bf V}\cdot{\bf N}\;\mbox{along }{\bf \Gamma}_{1}\quad\left(\mbox{\ensuremath{\therefore{\bf V}^{\perp}\cdot{\bf e}_{0}=0}}\right)
\end{align*}
(note that the decomposition still depends upon the boundary condition
for $\phi$ along ${\bf \Gamma}$). We choose $\delta{\bf W}={\bf V}^{\perp}$
to maintain the internal mass conservation constraint $\nabla\cdot{\bf W}=0$
as well as the flux constraints along ${\bf \Gamma}_{0}$ and ${\bf \Gamma}_{1}$.
This eliminates the dependencies upon $\lambda$ everywhere and upon
$\alpha$ along the temporal faces ${\bf \Gamma}_{0}$ and ${\bf \Gamma}_{1}$.
\begin{align*}
\delta E(\phi)= & \int_{{\bf \Omega}}\|{\bf V}^{\perp}\|^{2}+{\bf V}^{\perp}\cdot{\bf V}^{\parallel}-{\bf V}^{\perp}\cdot\nabla\lambda\;d{\bf X}+\int_{{\bf \Gamma}}\left(\lambda+\alpha\right)\,{\bf V}^{\perp}\cdot{\bf {\bf N}}+{\bf W}\cdot\left({\bf V}+\nabla_{S}\alpha\right)\delta{\bf \Gamma}\cdot{\bf N}\,dS\\
 & -\int_{{\bf \Gamma}_{0}}\left(\lambda+\alpha\right)\,\underbrace{{\bf V}^{\perp}\cdot{\bf e}_{0}}_{0}\,dS+\int_{{\bf \Gamma}_{1}}\left(\lambda+\alpha\right)\,\underbrace{{\bf V}^{\perp}\cdot{\bf e}_{0}}_{0}\,dS\\
= & \int_{{\bf \Omega}}\|{\bf V}^{\perp}\|^{2}+{\bf V}^{\perp}\cdot\nabla\phi+\lambda\underbrace{\nabla\cdot{\bf V}^{\perp}}_{0}\;d{\bf X}+\int_{{\bf \Gamma}}\alpha\,{\bf V}^{\perp}\cdot{\bf {\bf N}}+{\bf W}\cdot\left({\bf V}+\nabla_{S}\alpha\right)\delta{\bf \Gamma}\cdot{\bf N}\,dS\\
= & \int_{{\bf \Omega}}\|{\bf V}^{\perp}\|^{2}-\phi\underbrace{\nabla\cdot{\bf V}^{\perp}}_{0}\;d{\bf X}+\int_{{\bf \Gamma}}\left(\phi+\alpha\right)\,{\bf V}^{\perp}\cdot{\bf {\bf N}}+{\bf W}\cdot\left({\bf V}+\nabla_{S}\alpha\right)\delta{\bf \Gamma}\cdot{\bf N}\,dS\\
= & \int_{{\bf \Omega}}\|{\bf V}^{\perp}\|^{2}\;d{\bf X}+\int_{{\bf \Gamma}}\left(\phi+\alpha\right)\,{\bf V}^{\perp}\cdot{\bf {\bf N}}+{\bf W}\cdot\left({\bf V}+\nabla_{S}\alpha\right)\delta{\bf \Gamma}\cdot{\bf N}\,dS\\
= & \int_{{\bf \Omega}}\|{\bf V}-\nabla\phi\|^{2}\;d{\bf X}+\int_{{\bf \Gamma}}\left(\phi+\alpha\right)\,\left({\bf V}-\nabla\phi\right)\cdot{\bf {\bf N}}+{\bf W}\cdot\left({\bf V}+\nabla_{S}\alpha\right)\delta{\bf \Gamma}\cdot{\bf N}\,dS
\end{align*}

\subsection{Coupled boundary and flux perturbations\label{app:coupled}}

To maintain the vanishing flux constraint ${\bf W}\cdot{\bf N}=0$
along hypersurface ${\bf \Gamma}$ portion of the support boundary
$\partial{\bf \Omega}$, the normal perturbation $\delta{\bf \Gamma}\cdot{\bf N}$
of the boundary itself and the normal component of the solenoidal
field perturbation $\delta{\bf W}\cdot{\bf N}$ cannot be applied
independently but are coupled. This should not be surprising because
the field ${\bf W}$ defines the transport which, by virtue of determining
the density evolution, also determines the evolution of its support.
To determine the resulting coupling between a support perturbation
$\delta{\bf \Gamma}\cdot{\bf N}$ and the matching flux perturbation
$\delta{\bf W}\cdot{\bf N}$, we differentiate the following constraint
along the swept-out hypersurface
\begin{equation}
{\bf W}({\bf {\bf \Gamma}}(s))\cdot{\bf N}(s)=0\label{eq:noflux}
\end{equation}
where
\[
s=(s_{1},\ldots,s_{n})
\]
denotes isothermal coordinates aligned with the principal directions
${\bf T}_{1},\ldots,{\bf T}_{n}$ (unit tangent vectors) of the hypersurface.
We choose this parameterization for the convenient property that geodesic
torsions vanish along principal directions, and therefore
\[
\frac{\partial{\bf T}_{j}}{\partial s_{i}}=\begin{cases}
\kappa_{i}{\bf N}, & i=j\\
0, & i\ne j
\end{cases}
\]
were $\kappa_{i}$ denotes the principle curvature. We now expand
\[
\delta({\bf W}\cdot{\bf N})=\left(\frac{\partial{\bf W}}{\partial{\bf X}}\delta{\bf \Gamma}+\delta{\bf W}\right)\cdot{\bf N}+{\bf W}\cdot\delta{\bf N}=0
\]
to obtain
\begin{align*}
\delta{\bf W}\cdot{\bf N}= & -{\bf N}^{T}\frac{\partial{\bf W}}{\partial{\bf X}}\delta{\bf \Gamma}-{\bf W}\cdot\delta{\bf N}\\
= & \underbrace{-\left({\bf N}^{T}\frac{\partial{\bf W}}{\partial{\bf X}}{\bf N}\right)\delta{\bf \Gamma}\cdot{\bf N}-\sum_{i=1}^{n}\left({\bf N}^{T}\frac{\partial{\bf W}}{\partial{\bf X}}{\bf T}_{i}\right)\delta{\bf \Gamma}\cdot{\bf T}_{i}}_{{\bf \mbox{orthogonal components of }-{\bf N}^{T}\frac{\partial{\bf W}}{\partial{\bf X}}}\delta{\bf \Gamma}}-{\bf W}\cdot\underbrace{\left(\sum_{i=1}^{n}-{\bf T}_{i}\left(\frac{\partial(\delta{\bf \Gamma})}{\partial s_{i}}\cdot{\bf N}\right)\right)}_{\delta{\bf N}}\\
= & \left(\overbrace{\nabla\cdot{\bf W}}^{0}-{\bf N}^{T}\frac{\partial{\bf W}}{\partial{\bf X}}{\bf N}\right)\delta{\bf \Gamma}\cdot{\bf N}-\sum_{i=1}^{n}\left({\bf N}^{T}\frac{\partial{\bf W}}{\partial{\bf X}}{\bf T}_{i}\right)\delta{\bf \Gamma}\cdot{\bf T}_{i}+\sum_{i=1}^{n}{\bf W}\cdot{\bf T}_{i}\left(\frac{\partial(\delta{\bf \Gamma})}{\partial s_{i}}\cdot{\bf N}\right)\\
= & \left(\sum_{i=1}^{n}{\bf T}_{i}^{T}\frac{\partial{\bf W}}{\partial{\bf X}}{\bf T}_{i}\right)\delta{\bf \Gamma}\cdot{\bf N}-\sum_{i=1}^{n}\left({\bf N}^{T}\frac{\partial{\bf W}}{\partial{\bf X}}{\bf T}_{i}\right)\delta{\bf \Gamma}\cdot{\bf T}_{i}+\sum_{i=1}^{n}\left(\frac{\partial(\delta{\bf \Gamma})}{\partial s_{i}}\cdot{\bf N}\right)\left({\bf W}\cdot{\bf T}_{i}\right)\\
= & \left(\sum_{i=1}^{n}{\bf T}_{i}^{T}\overbrace{\frac{\partial{\bf W}}{\partial{\bf X}}{\bf T}_{i}}^{\frac{\partial{\bf W}}{\partial s_{i}}}+\overbrace{{\bf W}\cdot\frac{\partial{\bf T}_{i}}{\partial s_{i}}}^{\kappa_{i}{\bf W}\cdot{\bf N}=0}\right)\delta{\bf \Gamma}\cdot{\bf N}-\sum_{i=1}^{n}\left({\bf N}^{T}\frac{\partial{\bf W}}{\partial{\bf X}}{\bf T}_{i}\right)\delta{\bf \Gamma}\cdot{\bf T}_{i}+\sum_{i=1}^{n}\left(\frac{\partial(\delta{\bf \Gamma})}{\partial s_{i}}\cdot{\bf N}\right)\left({\bf W}\cdot{\bf T}_{i}\right)\\
= & \sum_{i=1}^{n}\frac{\partial}{\partial s_{i}}\left({\bf W}\cdot{\bf T}_{i}\right)\delta{\bf \Gamma}\cdot{\bf N}-\left({\bf N}^{T}\frac{\partial{\bf W}}{\partial{\bf X}}{\bf T}_{i}\right)\delta{\bf \Gamma}\cdot{\bf T}_{i}+\left(\frac{\partial(\delta{\bf \Gamma})}{\partial s_{i}}\cdot{\bf N}\right)\left({\bf W}\cdot{\bf T}_{i}\right)
\end{align*}
Differentiating the vanishing flux condition (\ref{eq:noflux}) along
each principal direction yields
\[
\frac{\partial}{\partial s_{i}}({\bf W}\cdot{\bf N})=\left(\frac{\partial{\bf W}}{\partial{\bf X}}{\bf T}_{i}\right)\cdot{\bf N}-{\bf W}\cdot\kappa_{i}{\bf T}=0
\]
allowing us to substitute
\[
{\bf N}^{T}\frac{\partial{\bf W}}{\partial{\bf X}}{\bf T}_{i}={\bf W}\cdot\kappa_{i}{\bf T}_{i}
\]
into our previous expression and to continue as follows
\begin{align*}
\delta{\bf W}\cdot{\bf N}= & \sum_{i=1}^{n}\frac{\partial}{\partial s_{i}}\left({\bf W}\cdot{\bf T}_{i}\right)\delta{\bf \Gamma}\cdot{\bf N}-\left({\bf W}\cdot\kappa_{i}{\bf T}_{i}\right)\delta{\bf \Gamma}\cdot{\bf T}_{i}+\left(\frac{\partial(\delta{\bf \Gamma})}{\partial s_{i}}\cdot{\bf N}\right)\left({\bf W}\cdot{\bf T}_{i}\right)\\
= & \sum_{i=1}^{n}\frac{\partial}{\partial s_{i}}\left({\bf W}\cdot{\bf T}_{i}\right)\delta{\bf \Gamma}\cdot{\bf N}-\left({\bf W}\cdot\kappa_{i}{\bf T}_{i}\right)\delta{\bf \Gamma}\cdot{\bf T}_{i}+\frac{\partial(\sum_{j=1}^{n}\left(\delta{\bf \Gamma}\cdot{\bf T}_{j}\right){\bf T}_{j}+\left(\delta{\bf \Gamma}\cdot{\bf N}\right){\bf N})}{\partial s_{i}}\cdot{\bf N}\left({\bf W}\cdot{\bf T}_{i}\right)\\
= & \sum_{i=1}^{n}\frac{\partial}{\partial s_{i}}\left({\bf W}\cdot{\bf T}_{i}\right)\delta{\bf \Gamma}\cdot{\bf N}-\kappa_{i}\left({\bf W}\cdot{\bf T}_{i}\right)\delta{\bf \Gamma}\cdot{\bf T}_{i}+\left(\sum_{j\ne i}\left(\delta{\bf \Gamma}\cdot{\bf T}_{j}\right)\underbrace{\frac{\partial{\bf T}_{j}}{\partial s_{i}}}_{0}+\left(\delta{\bf \Gamma}\cdot{\bf T}_{i}\right)\underbrace{\frac{\partial{\bf T}_{i}}{\partial s_{i}}}_{\kappa_{i}}+\frac{\partial}{\partial s_{i}}\left(\delta{\bf \Gamma}\cdot{\bf N}\right)\right)\left({\bf W}\cdot{\bf T}_{i}\right)\\
= & \sum_{i=1}^{n}\frac{\partial}{\partial s_{i}}\left({\bf W}\cdot{\bf T}_{i}\right)\delta{\bf \Gamma}\cdot{\bf N}+\frac{\partial}{\partial s_{i}}\left(\delta{\bf \Gamma}\cdot{\bf N}\right)\left({\bf W}\cdot{\bf T}_{i}\right)
\end{align*}
yielding the final expression for the coupling between $\delta{\bf W}\cdot{\bf N}$
and $\delta{\bf \Gamma}\cdot{\bf N}$.
\[
\delta{\bf W}\cdot{\bf N}=\sum_{i=1}^{n}\frac{\partial}{\partial s_{i}}\Bigl(({\bf W}\cdot{\bf T}_{i})(\delta{\bf \Gamma}\cdot{\bf N})\Bigr)
\]
It is generally difficult to invert this expression to express the
boundary perturbation $\delta{\bf \Gamma}\cdot{\bf N}$ as a function
of the flux perturbation $\delta{\bf W}\cdot{\bf N}$. However, in
the special case of zero flux perturbation ${\bf \delta{\bf W}\cdot{\bf N}=}0$,
we obtain
\[
\sum_{i=1}^{n}\frac{\partial}{\partial s_{i}}\Bigl(({\bf W}\cdot{\bf T}_{i})(\delta{\bf \Gamma}\cdot{\bf N})\Bigr)=0
\]
which, combined with the constraint that $\delta{\bf \Gamma}\cdot{\bf N}=0$
along the temporal boundaries of ${\bf \Gamma}$ at $t=0$ and at
$t=1$, only admits the solution $\delta{\bf \Gamma}\cdot{\bf N}=0$
along the entirety of the hypersurface ${\bf \Gamma}$. As such, a
vanishing flux perturbation implies a vanishing perturbation of the
support boundary. The converse is also trivially demonstrated directly
through equation (\ref{eq:noflux}).

\subsection{Gradient descent step size\label{sec:newton}}

\paragraph{Maximum allowed step factor}

Since any choice of step factor $\gamma_{k}$ in (\ref{eq:descent_step})
maintains the solenoidal constraint and boundary flux conditions for
${\bf W}_{k+1}$, we enforce the constraint $\rho_{k+1}(\text{{\bf X}})>0$
for all ${\bf X}\in{\bf \Omega}$ to determine the upper limit for
the allowable step factor $\gamma_{k}$.
\[
0<\rho_{k+1}={\bf W}_{k+1}\cdot{\bf e}_{0}=\left({\bf W}_{k}+\gamma_{k}\underbrace{\delta{\bf W}_{k}}_{-{\bf V}_{k}^{\perp}}\right)\cdot{\bf e}_{0}={\bf \rho}_{k}+\,\gamma_{k}\delta\rho_{k}
\]
Notice that this inequality is satisfied for arbitrarily large $\gamma_{k}$
whenever $\delta\rho_{k}\ge0$ and so an upper bound needs to be considered
only in cases where $\delta\rho_{k}({\bf X})<0$ for some ${\bf X}\in{\bf \Omega}$.
Accordingly, if we denote the set 
\[
{\bf \Omega}^{-}=\left\{ {\bf X}\in{\bf \Omega}\;|\;\delta\rho_{k}<0\right\} 
\]
we may formulate the following strict upper bound for $\gamma_{k}$.
\[
\gamma_{k,\max}=\begin{cases}
\min\limits _{{\bf \Omega}^{-}}\left(-{\displaystyle \frac{{\bf \rho}_{k}}{\delta\rho_{k}}}\right), & {\bf \Omega}^{-}\ne\emptyset\\
\infty, & {\bf \Omega}^{-}=\emptyset
\end{cases}
\]
A numerically robust way to compute $y_{k,\max}$ while avoiding numerical
overflow is to start with an exceedingly large estimate and then loop
through the points ${\bf X}\in{\bf \Omega}$, checking for the condition
$\gamma_{k,\max}\,\delta\rho_{k}<-\rho_{k}$. Whenever the condition
is detected, the value of $y_{k,\max}$ should be replaced by $-\frac{{\bf \rho}_{k}}{\delta\rho_{k}}$
at the detected location ${\bf X}$. Since this condition will never
be satisfied for points outside of the set ${\bf \Omega}^{+}$, there
is no need test whether ${\bf X}\in{\bf \Omega}^{-}$ beforehand.

\paragraph{Some preliminary calculations}

Using the dot notation for the derivative with respect to $\gamma_{k}$
we begin with a few preliminary calculations as follows. The last
two dot product expressions will be useful in formulating the Newton
update.

\begin{align*}
{\bf W}_{k+1} & ={\bf W}_{k}+\gamma_{k}\delta{\bf W}_{k}\\
\dot{{\bf W}}_{k+1} & =\delta{\bf W}_{k}\\
{\bf V}_{k+1} & ={\bf U}_{k+1}-\frac{1}{2}\left({\bf U}_{k+1}\cdot{\bf U}_{k+1}+1\right){\bf e}_{0}=\frac{{\bf W}_{k+1}}{{\bf W}_{k+1}\cdot{\bf e}_{0}}-\frac{1}{2}\left(\frac{{\bf W}_{k+1}\cdot{\bf W}_{k+1}}{\left({\bf W}_{k+1}\cdot{\bf e}_{0}\right)^{2}}+1\right){\bf e}_{0}\\
\dot{{\bf V}}_{k+1} & =\frac{1}{{\bf W}_{k+1}\cdot{\bf e}_{0}}\left(\delta{\bf W}_{k}-\frac{{\bf W}_{k+1}}{{\bf W}_{k+1}\cdot{\bf e}_{0}}\delta{\bf W}_{k}\cdot{\bf e}_{0}-\left(\frac{{\bf W}_{k+1}\cdot\delta{\bf W}_{k}}{{\bf W}_{k+1}\cdot{\bf e}_{0}}-\frac{{\bf W}_{k+1}\cdot{\bf W}_{k+1}}{\left({\bf W}_{k+1}\cdot{\bf e}_{0}\right)^{2}}\delta{\bf W}_{k}\cdot{\bf e}_{0}\right){\bf e}_{0}\right)\\
{\bf W}_{k+1}\cdot\dot{{\bf V}}_{k+1} & =0\\
\delta{\bf W}_{k}\cdot\dot{{\bf V}}_{k+1} & =\frac{1}{{\bf W}_{k+1}\cdot{\bf e}_{0}}\left(\delta{\bf W}_{k}\cdot\delta{\bf W}_{k}-2\left(\frac{\delta{\bf W}_{k}\cdot{\bf e}_{0}}{{\bf W}_{k+1}\cdot{\bf e}_{0}}\right){\bf W}_{k+1}\cdot\delta{\bf W}_{k}+\left(\frac{\delta{\bf W}_{k}\cdot{\bf e}_{0}}{{\bf W}_{k+1}\cdot{\bf e}_{0}}\right)^{2}{\bf W}_{k+1}\cdot{\bf W}_{k+1}\right)\\
 & =\frac{1}{{\bf W}_{k+1}\cdot{\bf e}_{0}}\left\Vert \delta{\bf W}_{k}-\frac{\delta{\bf W}_{k}\cdot{\bf e}_{0}}{{\bf W}_{k+1}\cdot{\bf e}_{0}}{\bf W}_{k+1}\right\Vert ^{2}=\frac{1}{{\bf W}_{k+1}\cdot{\bf e}_{0}}\left\Vert \underbrace{\delta{\bf W}_{k}-\left(\delta{\bf W}_{k}\cdot{\bf e}_{0}\right){\bf U}_{k+1}}_{{\bf e}_{0}\mbox{-components cancel}}\right\Vert ^{2}\\
 & =\frac{\left\Vert \delta{\rm p}_{k}-\delta\rho_{k}v_{k+1}\right\Vert ^{2}}{\rho_{k+1}}
\end{align*}
where $\delta{\rm p}_{k}$ (momentum increment) denotes the spatial
components of $\delta{\bf W}_{k}$ and where $\delta\rho_{k}$ (mass
increment) denotes its temporal component
\[
\delta{\bf W}_{k}=(\delta\rho_{k},\delta{\rm p}_{k})
\]

\paragraph{Optimal step factor}

We may employ Newton's method to maximize the descent taken by the
single step \ref{eq:descent_step}. To do so, we first derive expressions
for the first and second derivatives of the updated energy $E_{k+1}$
as a function of $\gamma_{k}$ (exploiting the preliminary calculations
listed above).

\begin{align*}
E_{k+1}(\gamma_{k}) & =\int_{{\bf \Omega}}{\bf W}_{k+1}(\gamma_{k})\cdot{\bf V}_{k+1}(\gamma_{k})\,d{\bf X}\\
\dot{E}_{k+1}(\gamma_{k}) & =\int_{{\bf \Omega}}{\bf \dot{W}}_{k+1}(\gamma_{k})\cdot{\bf V}_{k+1}(\gamma_{k})+\underbrace{{\bf W}_{k+1}(\gamma_{k})\cdot\dot{{\bf V}}_{k+1}(\gamma_{k})}_{0}\,d{\bf X}=\int_{{\bf \Omega}}\delta{\bf W}_{k}\cdot{\bf V}_{k+1}(\gamma_{k})\,d{\bf X}\\
\ddot{E}_{k+1}(\gamma_{k}) & =\int_{{\bf \Omega}}\delta{\bf W}_{k}\cdot\dot{{\bf V}}_{k+1}(\gamma_{k})\,d{\bf X}=\int_{{\bf \Omega}}\frac{\left\Vert \delta{\rm p}_{k}-\delta\rho_{k}v_{k+1}\right\Vert ^{2}}{\rho_{k+1}}\,d{\bf X}
\end{align*}
Notice that
\begin{align*}
\ddot{E}_{k+1}(\gamma_{k}) & \ge0\quad\forall\gamma_{k}\\
\dot{E}_{k+1}(0) & =\int_{{\bf \Omega}}\underbrace{\delta{\bf W}_{k}}_{-{\bf V}_{k}^{\perp}}\cdot{\bf V}_{k}\,d{\bf X}=\underbrace{-\int_{{\bf \Omega}}{\bf V}_{k}^{\perp}\cdot{\bf V}_{k}\,d{\bf X}=-\int_{{\bf \Omega}}\|{\bf V}_{k}^{\perp}\|^{2}\,d{\bf X}}_{\mbox{by orthogonality of \ensuremath{{\bf V}_{k}^{\perp}}and \ensuremath{{\bf V}_{k}^{\parallel}}}}\le0
\end{align*}
and so, by the first inequality, $E_{k+1}(\gamma_{k})$ is convex
with a unique local (\emph{i.e.} global) minimum which is, by the
second inequality, achieved for some $\gamma_{k}\ge0$. If, in addition,
\[
\dot{E}_{k+1}(\gamma_{k,\max})=\int_{{\bf \Omega}}\delta{\bf W}_{k}\cdot{\bf V}_{k+1}(\gamma_{k,\max})\,d{\bf X}\ge0
\]
then we can be apply Newton iterations $\gamma_{k}\to\gamma_{k}-\frac{\dot{E}_{k+1}}{\ddot{E}_{k+1}}$
to solve $\dot{E}_{k+1}(\gamma_{k})=0$ for $\gamma_{k}$ as follows.

\medskip{}

\begin{align*}
\mbox{Solve: } & \mbox{Poisson equation for \ensuremath{\lambda_{k}} given }{\bf V}_{k}\\
\mbox{Initialize: } & \delta{\bf W}_{k}={\bf V}_{k}^{\perp}={\bf V}_{k}-\nabla\lambda_{k}\\
 & \gamma_{k}=0\\
\mbox{Loop: } & {\bf W}_{k+1}={\bf W}_{k}+\gamma_{k}\delta{\bf W}_{k}=(\rho_{k+1},{\rm p}_{k+1})\\
 & {\bf V}_{k+1}=\left(-\frac{1}{2}\|v_{k+1}\|^{2}\,,\,v_{k+1}\right)\\
 & \Delta\gamma_{k}=\frac{\dot{E}_{k+1}}{\ddot{E}_{k+1}}=\frac{\int_{{\bf \Omega}}\delta{\bf W}_{k}\cdot{\bf V}_{k+1}\,d{\bf X}}{\int_{{\bf \Omega}}\frac{\left\Vert \delta\rho_{k}\,v_{k+1}-\delta{\rm p}_{k}\right\Vert ^{2}}{\rho_{k+1}}\,d{\bf X}}\\
 & \gamma_{k}\to\gamma_{k}-\Delta\gamma_{k}\;\left[\mbox{clip if needed so }0\le\gamma_{k}<\gamma_{k,\max}\right]
\end{align*}

\end{document}